\documentclass[12pt]{amsart}
\usepackage{amssymb}
\usepackage{amsfonts}
\usepackage{latexsym}
\usepackage{amscd}

\usepackage[all]{xy}

\input xy
\xyoption{all}

\numberwithin{equation}{section}

\newcommand{\Z}{{\mathbb{Z}}}
\newcommand{\N}{{\mathbb{N}}}

\newcommand{\R}{{\mathbb{R}}}
\newcommand{\pis}{unital purely infinite simple\ }

\newtheorem{lem}{Lemma}[section]
\newtheorem{corol}[lem]{Corollary}
\newtheorem{theor}[lem]{Theorem}
\newtheorem{prop}[lem]{Proposition}
\newtheorem{rema}[lem]{Remark}
\newtheorem{defi}[lem]{Definition}
\newtheorem{defis}[lem]{Definitions}

\newtheorem{p}[lem]{}

\def\p{^{\prime}}

\begin{document}
\title[Exchange Leavitt path algebras and stable rank]{Exchange Leavitt path algebras and stable rank}
\author{G. Aranda Pino}
\address{Departamento de \'Algebra, Geometr\'{\i}a y Topolog\'{\i}a,
Universidad de
M{\'a}laga,  29071 M{\'a}laga,
Spain.}\email{gonzalo@agt.cie.uma.es}
\author{E. Pardo}
\address{Departamento de Matem\'aticas, Universidad de C\'adiz,
Apartado 40, 11510 Puerto Real (C\'adiz),
Spain.}\email{enrique.pardo@uca.es}
\urladdr{http://www.uca.es/dept/matematicas/PPersonales/PardoEspino/index.HTML}
\author{M. Siles Molina}
\address{Departamento de \'Algebra, Geometr\'{\i}a y Topolog\'{\i}a,
Universidad de
M\'alaga,  29071 M\'alaga, Spain.}\email{mercedes@agt.cie.uma.es}
\thanks{The first author was partially supported by a FPU grant AP2001-1368
by the MEC. The first and third authors were partially supported
by the MCYT and Fondos FEDER, BFM2001-1938-C02-01,
MTM2004-06580-C02-02 and the ``Plan Andaluz de Investigaci\'on y
Desarrollo Tecnol\'ogico", FQM 336. The second author was
partially supported by the DGI and European Regional Development
Fund, jointly, through Project MTM2004-00149, by PAI III grant
FQM-298 of the Junta de Andaluc\'{\i}a, and by the Comissionat per
Universitats i Recerca de la Generalitat de Catalunya.}
\date{\today}
%
%
\begin{abstract}
We characterize Leavitt path algebras which are exchange rings in
terms of intrinsic properties of the graph and show that the
values of the stable rank for these algebras are $1$, $2$ or
$\infty$. Concrete criteria in terms of properties of the
underlying graph are given for each case.
\end{abstract}

\maketitle

\section*{Introduction}

For a row-finite graph $E$, the Leavitt path algebra $L(E)$ is the
algebraic analogue of the Cuntz-Krieger algebra $C^*(E)$ described
in \cite{R}. The pioneer papers where $L(E)$ is introduced and
studied are \cite{AMFP, AA1, AA2}. In \cite{AMFP}, Ara, Moreno and
Pardo carry out a study of the monoid $V(L(E))$. Concretely they
show that there is a natural isomorphism between the lattices of
graded ideals of $L(E)$ and that of order ideals of $V(L(E))$. In
\cite{AA1} and \cite{AA2} Abrams and Aranda Pino provide
characterizations of the simplicity and purely infinite
simplicity, respectively, of the Leavitt path algebra $L(E)$ in
terms of properties involving the graph $E$ only.

An associative unital ring $R$ is said to be an exchange ring if
$R_R$ has the exchange property introduced by Crawley and
J\'{o}nsson. The structure of exchange rings has been intensively
investigated by several authors; in the non necessarily unital
case, their study was initiated by Ara in \cite{Ara1}. On the
other hand, the concept of stable rank, introduced by Bass for
unital rings (see e.g. \cite{B}), is very useful in treating the
stabilization problem in K-theory. In \cite{Vas}, Vasershtein
opens with the definition of stable rank for a nonnecessarily
unital ring. For the more specific case of $C^*$-algebras, the
exchange property is closely related with the real rank: The
$C^*$-algebras having real rank zero are precisely those which are
exchange rings.

Following the philosophy of \cite{AMFP, AA1, AA2}, the aim of this
paper is to study the exchange property for Leavitt path algebras
and, focusing on this type of algebras, their stable rank.

Some of our sources of inspiration for the characterization of the
exchange property are the works of Jeong and Park \cite{JP} and
Bates, Hong, Raeburn and Szyma\'{n}ski \cite{BHRS}, while for the
stable rank it is the paper by Deicke, Hong and Szyma\'{n}ski
\cite{DHSz}. The proofs presented here significantly differ to
those of the analytic setting of $C^*$-algebras and the arguments
are necessarily different in the purely algebraic context since
many of the tools used there are not available in our case.



The paper is divided into seven sections. After some
preliminaries, we begin by stating basic properties concerning
special subsets of graphs. In particular, we study the ideals
generated by hereditary and saturated subsets of vertices and
cofinality of the graph.

Condition (K), studied in the third section, plays a central role
in the paper. On the one hand, it is precisely the condition we
need to impose on $E$ so that $L(E)$ is exchange; on the other
hand, the development of results concerning stable rank occur
under this hypothesis.

The main result characterizing exchange Leavitt path algebras
appears in Section 4:

\medskip

{\bf{Theorem \ref{tres condiciones}}.} {\it For a graph $E$, the
following conditions are equivalent:
\begin{enumerate}
\item $L(E)$ is an exchange ring.
\item $E/H$ satisfies condition (L) for every hereditary saturated subset
$H$ of $E^0$.
\item $E$ satisfies condition (K).
\item $\mathcal{L}_{\text{gr}}(L(E))=\mathcal{L}(L(E))$.
\item $E_H$ and $E/H$ satisfy condition (K) for every hereditary saturated
subset $H$ of $E^0$.
\item $E_H$ and $E/H$ satisfy condition (K) for some hereditary saturated
subset $H$ of $E^0$.
\end{enumerate}}

\medskip

The rest of the sections are devoted to compute the stable rank in
Leavitt path algebras satisfying condition (K). The first step
towards this aim is done in Section 5: First, by investigating the
absence of unital purely infinite simple quotients of $L(E)$
(Proposition \ref{noquot}). Secondly, by relating prime graded
ideals with maximal tails (Proposition \ref{gradedprime}). Then,
in Section 6, we calculate the stable rank for Leavitt path
algebras which do not have nonzero bounded graph traces and for
which every vertex lying on a closed simple path is left infinite
(Corollary \ref{elrizo}). The paper finishes in Section 7 with a
criterion to compute the stable rank for exchange Leavitt path
algebras:

\medskip

{\bf Theorem \ref{trychotomy}.} {\it Let $E$ be a graph satisfying
condition (K). Then, the values of the stable rank of $L(E)$ are:
\begin{enumerate}
\item $\mbox{sr}(L(E))=1$ if $E$ is acyclic. \item
$\mbox{sr}(L(E))=\infty$ if there exists $H\in \mathcal{H}_E$ such
that the quotient graph $E/H$ is nonempty, finite, cofinal and
contains no sinks. \item $\mbox{sr}(L(E))=2$ otherwise.
\end{enumerate}}


\section{Preliminaries}

Along this paper, we describe Leavitt path algebras following the
presentation of \cite[Sections 2 and 4]{AMFP} but using the
notation of \cite{AA1} for the elements.

A \emph{(directed) graph} $E=(E^0,E^1,r,s)$ consists of two countable sets
$E^0,E^1$ and maps $r,s:E^1 \to E^0$. The
elements of $E^0$ are called \emph{vertices} and the elements of $E^1$
\emph{edges}.

A vertex which emits no edges is called a  \emph{sink}. A graph
$E$ is \emph{finite} if $E^0$ is a finite set.  If $s^{-1}(v)$ is
a  finite set for every $v\in E^0$, then the graph is called
\emph{row-finite}. A \emph{path} $\mu$ in a graph $E$ is a
sequence of edges $\mu=(\mu_1, \dots, \mu_n)$ such that
$r(\mu_i)=s(\mu_{i+1})$ for $i=1,\dots,n-1$. In such a case,
$s(\mu):=s(\mu_1)$ is the \emph{source} of $\mu$ and
$r(\mu):=r(\mu_n)$ is the \emph{range} of $\mu$. An edge $e$ is an
exit for a path $\mu$ if there exists $i$ such that
$s(e)=s(\mu_i)$ and $e\neq \mu_i$. If $s(\mu)=r(\mu)$ and
$s(\mu_i)\neq s(\mu_j)$ for every $i\neq j$, then $\mu$ is a
called a \emph{cycle}. If $v=s(\mu)=r(\mu)$ and $s(\mu_i)\neq v$
for every $i>1$, then $\mu$ is a called a \emph{closed simple path
based at $v$}. We denote by $CSP_E(v)$ the set of closed simple
paths in $E$ based at $v$. For a path $\mu$ we denote by $\mu^0$
the set of its vertices, i.e., $\{s(\mu_1),r(\mu_i)\mid
i=1,\dots,n\}$. For $n\ge 2$ we define $E^n$ to be the set of
paths of length $n$, and $E^*=\bigcup_{n\ge 0} E^n$ the set of all
paths. We define a relation $\ge$ on $E^0$ by setting $v\ge w$ if
there is a path $\mu\in E^*$ with $s(\mu)=v$ and $r(\mu)=w$. A
subset $H$ of $E^0$ is called \emph{hereditary} if $v\ge w$ and
$v\in H$ imply $w\in H$. A hereditary set is \emph{saturated} if
every vertex which feeds into $H$ and only into $H$ is again in
$H$, that is, if $s^{-1}(v)\neq \emptyset$ and
$r(s^{-1}(v))\subseteq H$ imply $v\in H$. The set $T(v)=\{w\in
E^0\mid v\ge w\}$ is the \emph{tree} of $v$, and it is the
smallest hereditary subset of $E^0$ containing $v$. We extend this
definition for an arbitrary set $X\subseteq E^0$ by
$T(X)=\bigcup_{x\in X} T(x)$. Denote by $\mathcal{H}$ (or by
$\mathcal{H}_E$ when it is necessary to emphasize the dependence
on $E$) the set of hereditary saturated subsets of $E^0$. The
\emph{hereditary saturated closure} of a set $X$ is defined as the
smallest hereditary and saturated subset of $E^0$ containing $X$.
It is shown in \cite{AMFP} that the hereditary saturated closure
of a set $X$ is $\overline{X}=\bigcup_{n=0}^\infty \Lambda_n(X)$,
where
\begin{enumerate}
\item $\Lambda_0(X)=T(X)$, \item $\Lambda_n(X)=\{y\in E^0\mid
s^{-1}(y)\neq \emptyset$ and $r(s^{-1}(y))\subseteq
\Lambda_{n-1}(X)\}\cup \Lambda_{n-1}(X)$, for $n\ge 1$.
\end{enumerate}
Let $E=(E^0,E^1, r, s)$ be a  graph, and let $K$ be a field. We
define the {\em Leavitt path $K$-algebra} $L_K(E)$ associated with
$E$ ($L(E)$ when the based field is understood) as the $K$-algebra
generated by a set $\{v\mid v\in E^0\}$ of pairwise orthogonal
idempotents, together with a set of variables $\{e,e^*\mid e\in
E^1\}$, which satisfy the following relations:

(1) $s(e)e=er(e)=e$ for all $e\in E^1$.

(2) $r(e)e^*=e^*s(e)=e^*$ for all $e\in E^1$.

(3) $e^*e'=\delta _{e,e'}r(e)$ for all $e,e'\in E^1$.

(4) $v=\sum _{\{ e\in E^1\mid s(e)=v \}}ee^*$ for every $v\in E^0$
that emits edges.

Note that the relations above imply that $\{ee^*\mid e\in E^1\}$
is a set of pairwise orthogonal idempotents in $L(E)$. Note also
that if $E$ is a finite graph then we have $\sum _{v\in E^0} v=1$.
In general the algebra $L(E)$ is not unital, but it can be written
as a direct limit of unital Leavitt path algebras (with non-unital
transition maps), so that it is an algebra with local units. Along
this paper, we will be concerned only with row-finite graphs.

\section{Basic properties of graphs}

Let $E$ be a graph. For any subset $H$ of $E^0$, we will denote by $I(H)$
the ideal of $L(E)$ generated by $H$.

\begin{lem}\label{qualsevol}
If $H$ is a subset of $E^0$, then $I(H)=I({\overline{H}})$, and
$\overline{H}=I(H)\cap E^0$.
\end{lem}
\begin{proof}
Take $G=I(H)\cap E^0$. By \cite[Lemma 3.9]{AA1}, $G\in
\mathcal{H}$. Thus, by minimality, we get $H\subseteq
{\overline{H}}\subseteq G$, whence $I(H)\subseteq
I({\overline{H}})\subseteq I(G)$. Since $G\subseteq I(H)$, we have
$I(G)\subseteq I(H)$, so we get the desired equality. The second
statement holds by \cite[Proposition 4.2 and Theorem 4.3]{AMFP},
as desired.
\end{proof}

For a graph $E$ and a hereditary subset $H$ of $E^0$, we denote by $E/H$ the
\textit{quotient graph}
$$(E^0\setminus H, \{e\in E^1\mid r(e)\not\in H\}, r\vert_{(E/H)^1},
s\vert_{(E/H)^1}),$$ and by $E_H$ the \textit{restriction graph}
$$(H, \{e\in E^1\mid s(e)\in H\}, r\vert_{(E_H)^1},
s\vert_{(E_H)^1}).$$ Observe that while $L(E_H)$ can be seen as a
subalgebra of $L(E)$, the same cannot be said about $L(E/H)$.

Now, we recall that $L(E)$ has a $\mathbb{Z}$-grading. For every
$e\in E^1$, set the degree of $e$ as 1, the degree of $e^\ast$ as
-1, and the degree of every element in $E^0$ as 0. Then we obtain
a well-defined degree on the Leavitt path $K$-algebra $L(E)$,
thus, $L(E)$ is a $\mathbb{Z}$-graded algebra:
$$L(E)=\bigoplus\limits_{n\in \mathbb{Z}}L(E)_n, \quad L(E)_nL(E)_m\subseteq
L(E)_{n+m}, \  \hbox{for all }\ n, m \in \mathbb{Z}.$$

An ideal $I$ of a $\mathbb{Z}$-graded algebra $A=\bigoplus\limits_{n\in
\mathbb{Z}}A_n$ is a \emph{graded ideal} in
case $I=\bigoplus\limits_{n\in \mathbb{Z}} (I\cap A_n)$.

\begin{rema}\label{gr id}
{\rm $\mbox{ }$ An ideal $J$ of $L(E)$ is graded if and only if it
is generated by idempotents; in fact, $J=I(H)$, where $H=J\cap
E^0\in \mathcal{H}_E$. (See the proofs of \cite[Proposition 4.2
and Theorem 4.3]{AMFP}.)}
\end{rema}

\begin{lem}\label{quotient1}
Let $E$ be a graph and consider a proper $H\in \mathcal{H}_E$. Define $\Psi
: L(E)\rightarrow L(E/H)$ by setting $\Psi
(v)=\chi_{(E/H)^0}(v)v$, $\Psi (e)=\chi_{(E/H)^1}(e)e$ and $\Psi
(e^*)=\chi_{((E/H)^1)^*}(e^*)e^*$ for every vertex $v$
and every edge $e$, where $\chi_{(E/H)^0}:E^0\rightarrow K$ and
$\chi_{(E/H)^1}:E^1\rightarrow K$ denote the
characteristic functions. Then:
\begin{enumerate}
\item The map $\Psi$ extends to a $K$-algebra epimorphism of
$\mathbb{Z}$-graded algebras with $\mbox{Ker}(\Psi)=I(H)$ and
therefore $L(E)/I(H)\cong L(E/H)$. \item If $X$ is hereditary in
$E$, then $\Psi(X)\cap (E/H)^0$ is hereditary in $E/H$. \item For
$X\supseteq H$, $X\in \mathcal{H}_E$ if and only if $\Psi(X)\cap
(E/H)^0 \in \mathcal{H}_{(E/H)}$. \item For every $X\supseteq H$,
$\overline{\Psi(X)\cap (E/H)^0}= \Psi(\overline{X})\cap (E/H)^0$.

\end{enumerate}
\end{lem}
\begin{proof}
(1) It was shown in \cite[Proof of Theorem 3.11]{AA1} that $\Psi$
extends to a $K$-algebra morphism. Since $H\in \mathcal{H}_E$,
$\Psi $ extends to a well-defined morphism. By definition, $\Psi$
is $\Z$-graded and onto. Moreover, $I(H)\subseteq \mbox{Ker}(\Psi
)$.

Since $\Psi $ is a graded morphism, $\mbox{Ker}(\Psi )\in
\mathcal{L}_{\text{gr}}(L(E))$. By \cite[Theorem 4.3]{AMFP}, there
exists $X\in \mathcal{H}_E$ such that $\mbox{Ker}(\Psi )=I(X)$. By
Lemma \ref{qualsevol}, $H= I(H)\cap E^0\subseteq I(X)\cap E^0=X$.
Hence, $I(H)\ne \mbox{Ker}(\Psi )$ if and only if there exists
$v\in X\setminus H$. But then $\Psi (v)=v\ne 0$ and $v\in
\mbox{Ker}(\Psi )$, which is impossible.

(2) It is clear by the definition of $\Psi$.

(3) Since $\Psi $ is a graded epimorphism, there is a bijection
between graded ideals of $L(E/H)$ and graded ideals of $L(E)$
containing $I(H)$. Thus, the result holds by \cite[Theorem
4.3]{AMFP}.

(4) It is immediate by part (3).
\end{proof}

Recall that a ring $R$ is said to be an \emph{idempotent ring} if
$R=R^2$. For an idempotent ring $ R $ we denote by $R-$Mod the
full subcategory of the category of all left $R$-modules whose
objects are the ``unital" nondegenerate modules. Here a left
$R$-module $M$ is said to be \emph{unital} if $M=RM,$ and $M$ is
said to be \emph{nondegenerate} if, for $m\in M$, $Rm=0$ implies
$m=0$. Note that if $R$ has an identity then $R-$Mod is the usual
category of left $R-$modules.

We will use the well-known definition of a Morita context in the case where
the rings $ R $ and $ S$ have not
necessarily an identity. Let $R$ and $S$ be idempotent rings. We say that
$(R, S,M, N,\varphi,\psi ) $ is a
\emph{(surjective) Morita context} if $ _RM_S $ and $ _SN_R $ are unital
bimodules and $\varphi : N{\otimes}_R M \to S,
$ $\psi : M{\otimes}_S N\to R$ are surjective $S$-bimodule and $R$-bimodule
maps, respectively, satisfying the
compatibility relations: $\varphi (n\otimes m) n\p =n\psi (m\otimes  n\p ),
$ $m\p \varphi (n\otimes m)= \psi (m\p
\otimes n)m $ for every $ m, m\p \in M,$ $ n, n\p\in N. $

In \cite{GS} (see Proposition 2.5 and Theorem 2.7) it is proved that if  $ R
$ and $ S$ are two idempotent rings, then
$R-$Mod and $S-$Mod are equivalent categories if and only if there exists a
(surjective) Morita context $ (R, S, M,
N,\varphi ,\psi ) $. In this case, we will say that the rings $R$ and $S$
are \emph{Morita equivalent} and we will
refer to as the (surjective) Morita context  $(R, S,M, N)$.

\begin{lem}\label{Morita}
Let $E$ be a graph and $H\subseteq E^0$ a proper hereditary subset. Then
$L(E_H)$ is Morita equivalent to $I(H)$.
\end{lem}
\begin{proof}
Let $ H=\{v_i\mid i\geq 1\}$, and consider the ascending family of
idempotents $e_n=\sum_{i=1}^{n}v_i$, ($n\geq 1$). By \cite[Lemma
1.6]{AA1}, $\{e_n\mid n\geq 1\}$ is a set of local units for
$L(E_H)$, so that $L(E_H)=\bigcup_{i\geq 1}e_i L(E)e_i$. Since
$I(H)$ is generated by the idempotents $v_i\in H$, it is a
non-degenerated idempotent ring. Moreover, $I(H)=\bigcup_{i\geq
1}L(E)e_iL(E)$. It is not difficult to see that $(\sum_{i\geq
1}e_iL(E)e_i, \sum_{i\geq 1}L(E)e_iL(E), \sum_{i\geq 1}L(E)e_i,
\sum_{i\geq 1}e_iL(E))$ is a (surjective) Morita context for the
idempotent rings $\sum_{i\geq 1}e_iL(E)e_i=L(E_H)$ and
$\sum_{i\geq 1}L(E)e_iL(E)=I(H)$, hence $I(H)$ is Morita
equivalent to $L(E_H)$.
\end{proof}

Under certain conditions we will see in Section 5 that $I(H)$ is
not only Morita equivalent to a Leavitt path algebra; in fact it
is isomorphic to a Leavitt path algebra.

\begin{lem}\label{projsather}
Let $H\in \mathcal{H}_E$, and let $X\subseteq H$ be any subset. Then, $X\in
\mathcal{H}_E$ if and only if $X\in
\mathcal{H}_{E_H}$.
\end{lem}
\begin{proof}
First, suppose that $X\in \mathcal{H}_{E_H}$. Since $X\subseteq H$, if $e\in
E^1$ and $s(e)\in X\subseteq H$, we have
$e\in E_H^1$. Hence, $r(e)\in X$, so that $X$ is hereditary into $E$. Now,
let $v\in E^0$ such that $s_E^{-1}(v)\ne
\emptyset$ and $r_E(s_E^{-1}(v))\subseteq X\subseteq H$. Then, $v\in H$, so
that $s_{E_H}^{-1}(v)=s_E^{-1}(v)$ (in
particular, it is nonempty), and
$r_{E_H}(s_{E_H}^{-1}(v))=r_E(s_E^{-1}(v))\subseteq X$, whence $v\in X$.
Thus, $X\in
\mathcal{H}_E$.

Suppose that $X\in \mathcal{H}_E$. Clearly, $X$ is hereditary in
$E_H$ (because $X\subseteq H$). Now, let $v\in H$ such that
$s_{E_H}^{-1}(v)\ne \emptyset$ and
$r_{E_H}(s_{E_H}^{-1}(v))\subseteq X\subseteq H$. Since $v\in H$,
$s_E^{-1}(v)=s_{E_H}^{-1}(v)$ (in particular, it is nonempty), and
$r_E(s_E^{-1}(v))=r_{E_H}(s_{E_H}^{-1}(v))\subseteq X$, whence
$v\in X$. Thus, $X\in \mathcal{H}_{E_H}$.
\end{proof}

\begin{lem}\label{gen_epi}
Let $E$ be a graph and $H\in \mathcal{H}_E$. Then, the canonical map
$$K_0(L(E))\rightarrow K_0(L(E)/I(H))$$ is an
epimorphism.
\end{lem}
\begin{proof}
If $H=E^0$ or $H=\emptyset$, the result follows trivially. Now, suppose $H$
to be a proper subset of $E^0$. By Lemma
\ref{quotient1} (1) we have $L(E)/I(H)\cong L(E/H)$. By \cite[Lemma
5.6]{AMFP},
$$V(L(E))/V(I(H))\cong V(L(E/H))\cong V(L(E)/I(H)).$$
Since $L(E)$ and $L(E/H)$ have a countable unit, we have that
$K_0(L(E))=\mbox{Grot}(V(L(E)))$ and
$K_0(L(E/H))=\mbox{Grot}(V(L(E/H)))$. Hence, the canonical map
$K_0(L(E))\rightarrow K_0(L(E)/I(H))$ is clearly an
epimorphism, as desired.
\end{proof}

We denote by $E^\infty$ the set of infinite paths
$\gamma=(\gamma_n)_{n=1}^\infty$ of the graph $E$ and by $E^{\le
\infty}$ the set $E^\infty$ together with the set of finite paths
in $E$ whose end vertex is a sink. We say that a vertex $v$ in a
graph $E$ is \emph{cofinal} if for every $\gamma\in E^{\le
\infty}$ there is a vertex $w$ in the path $\gamma$ such that
$v\ge w$. We say that a  graph $E$ is \emph{cofinal} if so are all
the vertices of $E$.

Observe that if a graph $E$ has cycles, then $E$ cofinal implies that every
vertex connects to a cycle.

\begin{lem}\label{ojito}
If $E$ is cofinal, and $v\in E^0$ is a sink, then:
\begin{enumerate}
\item The only sink of $E$ is $v$.
\item For every $w\in E^0$, $v\in T(w)$.
\item $E$ contains no infinite paths. In particular, $E$ is acyclic.
\end{enumerate}
\end{lem}
\begin{proof}\mbox{ }
\begin{enumerate}
\item It is obvious from the definition.
\item Since $T(v)=\{ v\}$, the result follows from the definition of $T(v)$
by
considering the path $\gamma=v\in E^{\leq \infty}$.
\item If $\alpha \in E^{\infty }$, then there exists $w\in \alpha
^0$ such that $v\geq w$, which is impossible. Thus, in particular, $E$
contains no closed simple paths, and therefore
no cycles.
\end{enumerate}
\end{proof}

Next result is known in the case of graphs without sinks. Since we
have no knowledge of the existence of a (published) version of the
result in the general case, we give a proof for the sake of
completeness.

\begin{lem}\label{cofinohersat}
A graph $E$ is cofinal if and only if $\mathcal{H}=\{ \emptyset , E^0 \}$.
\end{lem}
\begin{proof}
Suppose $E$ to be cofinal. Let $H\in \mathcal{H}$ with $\emptyset
\ne H \ne E^0$. Fix $v\in E^0\setminus H$ and build a path
$\gamma\in E^{\le \infty}$ such that $\gamma^0\cap H=\emptyset$:
If $v$ is a sink, take $\gamma=v$. If not, then $s^{-1}(v) \neq
\emptyset$ and $r(s^{-1}(v))\nsubseteq H$; otherwise, $H$
saturated implies $v\in H$, which is impossible. Hence, there
exists $e_1\in s^{-1}(v)$ such that $r(e_1)\notin H$. Let
$\gamma_1=e_1$ and repeat this process with $r(e_1)\not\in H$. By
recurrence either we reach a sink or we have an infinite path
$\gamma$ whose vertices are not in $H$, as desired. Now consider
$w\in H$. By the hypothesis, there exists $z\in\gamma$ such that
$w\ge z$, and by hereditariness of $H$ we get $z\in H$,
contradicting the definition of $\gamma$.

Conversely, suppose that $\mathcal{H}=\{ \emptyset , E^0 \}$. Take
$v\in E^0$ and $\gamma \in E^{\le \infty}$, with $v\not\in
\gamma^0$ (the case $v\in \gamma^0$ is obvious). By hypothesis the
hereditary saturated subset generated by $v$ is $E^0$, i.e.,
$E^0=\bigcup_{n\ge 0} \Lambda_n(v)$. Consider $m$, the minimum $n$
such that $\Lambda_n(v)\cap \gamma^0 \ne \emptyset$, and let $w
\in\Lambda_m(v)\cap \gamma^0$. If $m>0$, then by minimality of $m$
it must be $s^{-1}(w)\ne \emptyset$ and $r(s^{-1}(w))\subseteq
\Lambda_{m-1}(v)$. The first condition implies that $w$ is not a
sink and since $\gamma=(\gamma_n)\in E^{\le \infty}$, there exists
$i\geq 1$ such that $s(\gamma_i)=w$ and $r(\gamma_i)=w'\in
\gamma^0$, the latter meaning that $w'\in r(s^{-1}(w))\subseteq
\Lambda_{m-1}(v)$, contradicting the minimality of $m$. Therefore
$m=0$ and then $w\in \Lambda_0(v)=T(v)$, as we needed.
\end{proof}

\section{Condition (K)}

We begin this section by recalling the two following well-known
notions which will play a central role in the sequel.
\begin{enumerate}
\item A graph $E$ satisfies condition (L) if every closed simple path has an
exit, equivalently \cite[Lemma 2.5]{AA1},
if every cycle has an exit.
\item A graph $E$ satisfies condition (K) if for each vertex $v$ on a closed
simple path there exists at least two distinct closed simple paths
$\alpha, \beta$ based at $v$, or, following \cite{AA2},
$V_1=\emptyset$.
\end{enumerate}

\begin{rema}\label{nova2}
{\rm $\mbox{ }$
\begin{enumerate}
\item Notice that if $E$ satisfies condition (K) then it satisfies
condition (L).
\item According to \cite[Lemma 7]{AA2}, if $L(E)$ is simple then it
satisfies condition (K).
\end{enumerate}}
\end{rema}

It is not difficult to see that if $E$ satisfies condition (L)
then so does $E_H$, whereas $E/H$ need not. Condition (K) has a
better behaviour as it is shown in the following result.

\begin{lem}\label{sub_quot_graph}
Let $E$ be a graph and $H$ a hereditary subset of $E^0$. If $E$
satisfies condition (K), so do $E_H$ and $E/H$.
\end{lem}
\begin{proof}
We will see $CSP_E(v)=CSP_{E_H}(v)$ and $CSP_E(w)=CSP_{E/H}(w)$
for every $v\in H$ and $w\in E^0\setminus H$. Clearly,
$CSP_{E_H}(v)\subseteq CSP_E(v)$; conversely, let $\alpha \in
CSP_E(v)$, and suppose $\alpha =(\alpha_1, \dots ,\alpha_n)$.
Since $H$ is hereditary and $s(\alpha_1)=v\in H$, we get
$r(\alpha_1)=s(\alpha_2)\in H$. Thus, by recurrence, $\alpha\in
CSP_{E_H}(v)$ and the result holds.

Now, let $v\in E^0\setminus H$ and consider $\alpha =(\alpha_1,
\dots ,\alpha_n)\in CSP_E(v)$. Since $r(\alpha _n)=v\not\in H$ we
get $\alpha _n\in (E/H)^1$.  If $\alpha_{n-1}\notin (E/H)^1$ then
$r(\alpha_{n-1})=s({\alpha_n})\in H$ and $H$ hereditary implies
$v=r({\alpha_n})\in H$, a contradiction. By recurrence, $\alpha\in
CSP_{E/H}(v)$; since the converse is immediate, the result
follows.
\end{proof}

For an algebra $A$, denote by $\mathcal{L}(A)$ and
$\mathcal{L}_{\text{gr}}(A)$ the lattices of ideals and graded
ideals, respectively, of $A$. The following proposition provides a
description of the ideals of $L(E)$ for $E$ a graph
satisfying condition (K).

\begin{prop}\label{iso_ideal}
If a graph $E$ satisfies condition (K) then, for every ideal $J$
of $L(E)$, $J=I(H)$, where $H=J\cap E^0$ is a hereditary saturated
subset of $E^0$. In particular, $\mathcal{L}_{\text{gr}}(L(E))=
\mathcal{L}(L(E))$.
\end{prop}
\begin{proof}
Let $J$ be a nonzero ideal of $L(E)$. By \cite[Lemma 3.9]{AA1}
(which can be applied because $E$ satisfies condition (L) by
Remark \ref{nova2} (1)) and \cite[Proposition 6]{AA2}, $H=J\cap
E^0\neq\emptyset$ is a hereditary saturated subset of $E^0$.
Therefore, and taking into account Remark \ref{gr id}, $I(H)$ is a
graded ideal of $L(E)$ contained in J.

Suppose $I(H)\ne J$. Then, by Lemma \ref{quotient1} (1), $$0\ne
J/I(H)\lhd L(E)/I(H)\cong L(E/H).$$ Thus, $E/H$ satisfies
condition (L) by Lemma \ref{sub_quot_graph} and Remark \ref{nova2}
(1). Now, consider the isomorphism (of $K$-algebras)
$\overline{\Psi}: L(E)/I(H) \to L(E/H)$ given by
$\overline{\Psi}(x+I(H))=\Psi(x)$ (for $\Psi$ as in Lemma
\ref{quotient1}). By \cite[Proposition 6]{AA2}, $\emptyset \neq
\overline{\Psi}(J/I(H))\cap (E^0\setminus H)=\Psi(J) \cap
(E^0\setminus H)$, so there exists $v\in J$ with $\Psi(v)\in
\Psi(J)$. But $v\in E^0 \cap J=H$ and, on the other hand,
$v=\Psi(v)\in E^0\setminus H,$ which is impossible.

To finish, take into account that $J$ is an ideal generated by idempotents
and apply Remark \ref{gr id}.
\end{proof}

In the following section the converse of the previous result is
proved.

\begin{corol}\label{K-epi}
If $E$ satisfies condition (K), then for every ideal $I$ of $L(E)$, the
canonical map
$$K_0(L(E))\rightarrow K_0(L(E)/I)$$ is an epimorphism.
\end{corol}
\begin{proof}
By Proposition \ref{iso_ideal}, $I=I(H)$ for the hereditary saturated subset
$H=I\cap E^0$ of $E^0$. Then, the result
holds by Lemma \ref{gen_epi}.
\end{proof}

Recall that a \emph{matricial algebra} is a finite direct product
of full matrix algebras over $K$, while a \emph{locally matricial
algebra} is  a direct limit of matricial algebras.

The following result can be obtained as a corollary of Proposition
\ref{iso_ideal}. However we do not include its proof because it
can be reached by doing slight changes on that of \cite[Corollary
2.3]{KPR}.

\begin{corol}\label{matricial}
If $E$ is a finite acyclic graph, then $L(E)$ is a $K$-matricial algebra.
\end{corol}

\begin{corol}\label{ultramatricial}
If $E$ is an acyclic graph, then $L(E)$ is a locally matricial
$K$-algebra.
\end{corol}
\begin{proof}
By \cite[Lemma 2.2]{AMFP}, $L(E)\cong \varinjlim L(X_n)$, where
$X_n$ is a finite subgraph of $E$ for all $n\geq 1$. Hence, $X_n$
is a finite acyclic graph for every $n\geq 1$, whence the result
holds by Corollary \ref{matricial}.
\end{proof}

Recall that a \emph{graph homomorphism} $f\colon E=(E^0,E^1)\to
F=(F^0,F^1)$ is given by two maps $f^0\colon E^0\to F^0$ and
$f^1\colon E^1\to F^1$ such that $r_F(f^1(e))=f^0(r_E(e))$ and
$s_F(f^1(e))=f^0(s_E(e))$ for every $e\in E^1$. We say that a
graph homomorphism $f$ is {\em complete} in case $f^0$ is
injective and $f^1$ restricts to a bijection from $s_E^{-1}(v)$
onto $s_F^{-1}(f^0(v))$ for every $v\in E^0$ that emits edges.
Note that under the assumptions above, the map $f^1$ must also be
injective.

\begin{lem}\label{finalkey}
If $E$ is a graph satisfying condition (K) then there exists an
ascending family $\{X_n\}_{n\geq 0}$ of finite subgraphs such
that:
\begin{enumerate}
\item For every $n\geq 0$, $X_n$ satisfies condition (K).
\item For every $n\geq 0$, the inclusion map $X_n\subseteq E$ is a complete
graph homomorphism.
\item $E=\bigcup_{n\geq 0}X_n$.
\end{enumerate}
\end{lem}
\begin{proof}
We will construct $X_n$ by recurrence on $n$. First, we enumerate
$E^0=\{v_n\mid n\geq 0\}$. Then, we define $X_0=\{v_0\}$. Clearly,
$X_0$ satisfies condition (K) and also $X_0\subseteq E$ is a
complete graph homomorphism.

Now, suppose we have constructed $X_0, X_1, \dots , X_n$ satisfying $(1)$
and $(2)$. Consider the graph ${\widetilde
X}_{n+1}$ with: (a) ${\widetilde X}_{n+1}^1=X_n^1\cup \{e\in E^1\mid s(e)\in
X_n^0\}$; (b) ${\widetilde
X}_{n+1}^0=X_n^0\cup \{v_{n+1}\}\cup \{r(e)\mid e\in {\widetilde
X}_{n+1}^1\}$. Clearly, ${\widetilde X}_{n+1}$ is
finite and satisfies $(2)$. If it also satisfies $(1)$, we define
$X_{n+1}={\widetilde X}_{n+1}$.

Suppose that ${\widetilde X}_{n+1}$ does not satisfy condition
(K). Consider the set of all cycles based at vertices in
${\widetilde X}_{n+1}$, $\mu_1^1, \dots ,\mu_1^k\subseteq
{\widetilde X}_{n+1}$ such that: (i) $\mu_1^i\nsubseteq X_n$ for
any $1\leq i \leq k$; (ii) for every $1\leq i \leq k$ and every
$v\in \mu_1^i$, $\mbox{card}(CSP_{{\widetilde X}_{n+1}}(v))=1$.
Since ${\widetilde X}_{n+1}\subseteq E$ and $E$ satisfies
condition (K), there exist closed simple paths $\mu_2^1,\dots
,\mu_2^k\subseteq E$ such that, for each $1\leq i \leq k$,
$\mu_1^i\neq \mu_2^i$ and $\mu_1^i\cap \mu_2^i\ne \emptyset$. For
each $1\leq i \leq k$, let $\mu_2^i=(e_1^i, \dots , e_{j_i}^i)$.

We consider the finite subgraph ${\widetilde Y}_{n+1}$ of $E$ such that: (a)
${\widetilde Y}_{n+1}^1={\widetilde
X}_{n+1}^1\cup \{e_l^i\mid 1\leq i \leq k, 1\leq l \leq j_i\}$; (b)
${\widetilde Y}_{n+1}^0={\widetilde X}_{n+1}^0\cup
\{s(e_l^i),r(e_l^i)\mid 1\leq i \leq k, 1\leq l \leq j_i\}$. Clearly,
${\widetilde Y}_{n+1}$ satisfies $(1)$.

Now, let $X_{n+1}$ be the finite subgraph of $E$ such that: (a)
$X_{n+1}^1= {\widetilde Y}_{n+1}^1\cup \{f\in E^1\mid s(f)\in
(\mu_2^i)^0 \mbox{ for some } 1\leq i\leq k\}$; (b)
$X_{n+1}^0={\widetilde Y}_{n+1}^0\cup \{ r(e)\mid e\in
X_{n+1}^1\}$. If $\mu\subseteq X_{n+1}$ is a closed simple path
such that $\mu\nsubseteq {\widetilde Y}_{n+1}$, then either it
appears because one of the $e\in X_{n+1}^1\setminus {\widetilde
Y}_{n+1}^1$ is a single loop based at some vertex in one
$\mu_2^i$, or $s(e)\in (\mu_2^i)^0$ and $r(e)$ connects to a path
that comes back to $s(e)$. In any case, the (potential) new closed
simple paths are based at vertices of $\mu_2^i$ for some $i$,
whence $X_{n+1}$ satisfies $(1)$. Also, since the step from
${\widetilde Y}_{n+1}$ to $X_{n+1}$ adds all the exits of all the
vertices in the cycles $\mu_2^i$, we conclude that for any vertex
$v\in X_{n+1}^0$, $v$ is either a sink, or every $e\in E^1$ with
$s(e)\in X_{n+1}^0$ belongs to $X_{n+1}^1$. Hence,
$X_{n+1}\subseteq E$ is a complete graph homomorphism. This
completes the recurrence argument.

Finally, since $v_n\in X_n$ for every $n\geq 0$, we conclude that
$E^0=\bigcup_{n\geq 0}X_n^0$ and by $(2)$, $E^1=\bigcup_{n\geq
0}X_n^1$.
\end{proof}

The following definitions can be found in \cite[Definition 3.2]{JP}. Let $F$
be a subgraph of a graph $E$. Then:
\begin{enumerate}
\item The \emph{loop completion} $\ell _{E}(F)$  of $F$ in $E$ is the
subgraph of $E$ obtained as the union of $F$
with every closed path based at an element of $F^0$. \item The \emph{exit
completion} $F_e$ of $F$ is a subgraph
obtained by adding to $F$ the edges $V=\{e\in E^1\mid s(e)=s(f) \mbox{ for
some } f\in F^1\}$, and the vertices $\{
r(e)\mid e\in V\}$. We say that $F$ is \emph{exit complete} if $F=F_e$.
\end{enumerate}

\begin{lem}\label{medkey}
If $F$ is an exit complete subgraph of a graph $E$, then $L(F)$ is
isomorphic to a subalgebra of $L(E)$.
\end{lem}
\begin{proof}
Denote by $e'$ the edges of $F$ seen inside $E$,  and by $v'$ the
vertices of $F$ seen inside $E$. Since $F=F_e$, for every vertex
$v\in F^0$ we have that either $v$ is a sink or
$s_F^{-1}(v)=s_E^{-1}(v')$. Thus, the relations defining $L(F)$
and $L(F')\subseteq L(E)$ are exactly  the same, so that there is
a natural injective morphism form $L(F)$ to $L(E)$, as desired.
\end{proof}

\begin{lem}\label{finitecompletions}
If $F$ is a subgraph of a graph $E$ then:
\begin{enumerate}
\item \label{cond1}$F_e$ is exit complete.
\item \label{cond2}If $F$ is finite then so is $F_e$ whereas $l_E(F)$
need not be.
\end{enumerate}
\end{lem}

\begin{proof}$\mbox{ }$

(1) Clearly $F^1_e\subseteq (F_e)^1_e$. Let us see the other
inclusion. Take $g \in (F_e)^1_e$. If $g\in F^1_e$ we have
finished. If not, there exists $f\in F^1_e$ with $s(g)=s(f)$. We
have two possibilities: If $f\in F^1$ then, by definition, $g\in
F^1_e$. If $f\not\in F^1$ we can find $h\in F^1$ for which
$s(f)=s(h)$. Therefore $s(g)=s(h)$ and again $g\in F^1_e$. Now it
easily follows $F^0_e=(F_e)^0_e$. \par (2) Since $F$ is finite
(and row-finite) then $F^1$ is finite. Now, for each $f\in F^1$
there are finitely many edges $e\in E^1$ with $s(e)=s(f)$ (because
$E$ is row-finite), and therefore we are adding a finite number of
edges and consequently of vertices. Thus, $F_e$ is finite. To show
that $l_E(F)$ can be infinite, consider the infinite graph $E$
$$\xymatrix{{\bullet}^v \ar@/^4pt/ [r] & {\bullet} \ar@/^4pt/ [r] \ar@/^4pt/
[l] & {\bullet} \ar@/^4pt/ [r] \ar@/^4pt/ [l]
& {\bullet} \ar@/^4pt/@{.} [r] \ar@/^4pt/ [l] & {\ldots} \ar@/^4pt/@{.}
[l]}$$
Then $F=(\{v\},\emptyset)$ is finite while $l_E(F)=E$ is not.
\end{proof}

\begin{lem}\label{conditionKcompletions}
Let $E$ be a graph and $T$ be any subgraph of $E$. Define $F=l_E(T)$,
$G=l_E(T)_e$, $S$ the set of sinks of $G$ and
$J=G/\overline{S}$. Then:
\begin{enumerate}
\item \label{grafoF} $CSP_F(v)=CSP_E(v)$ for every $v\in F^0$.
\item \label{grafoG} $CSP_G(v)=CSP_E(v)$ for every $v\in G^0$ such
that $CSP_G(v)\neq \emptyset$.
\item \label{grafoJ} $CSP_J(v)=CSP_E(v)$ for every $v\in J^0$ such that
$CSP_J(v)\neq
\emptyset$.
\item \label{condicionK} If $E$ satisfies condition (K) then so do $F$, $G$
and $J$.
\end{enumerate}
\end{lem}
\begin{proof} (\ref{grafoF}) is evident from the definition of the loop
completion.

(\ref{grafoG}). Consider $v\in G^0$ such that $CSP_G(v)\neq
\emptyset$ and take $p=(p_1,\dots,p_n)\in CSP_G(v)$. Suppose
$p^0\cap F^0= \emptyset$; then for an arbitrary edge $p_i$ we have
$r(p_i)\not\in F^0$. The construction of the exit completion
yields that $p_i$ is a new added edge and consequently it exists
$f\in F^1$ with $s(f)=s(p_i)$ and hence $s(p_i)\in p^0\cap F^0$, a
contradiction. Therefore $p^0\cap F^0\neq \emptyset$. Take $w$ in
the previous intersection. Then $CSP_G(v)=$ ($v$ and $w$ are in
the same closed path) $CSP_G(w)\supseteq $ (because $F$ is a
subgraph of $G$) $CSP_{F}(w)=$ (by (\ref{grafoF})) $CSP_E(w)=$
($v$ and $w$ are in the same closed path) $CSP_E(v)\supseteq
CSP_G(v)$. Hence, (\ref{grafoG}) holds. Note that the result may
fail for $CSP_G(v)=\emptyset$.

(\ref{grafoJ}) Let $v\in J^0$ such that $CSP_J(v)\neq \emptyset$.
Obviously $CSP_J(v)\subseteq CSP_G(v)$. Now consider
$p=(p_1,\dots,p_k)\in CSP_G(v)$. We claim that $p^0\cap
\overline{S}=\emptyset$. If not, there exists $m=\min\{n\in
\mathbb{N}\mid p^0\cap \Lambda_n(S)\neq \emptyset\}$. Take $v\in
p^0\cap \Lambda_m(S)$. If $m>0$ then by minimality we have that
$r(s^{-1}(v))\subseteq \Lambda_{m-1}(S)$. In particular, if
$v=s(p_i)$ then $r(p_i)\in p^0\cap \Lambda_{m-1}(S)$, which
contradicts the minimality of $m$. If $m=0$ then $v\in p^0\cap
\Lambda_0(S)=p^0\cap S=$ ($S$ is the set of sinks) $p^0\cap S$.
This is absurd since $p$ has no sinks. Any possibility leads to a
contradiction so $p^0\subseteq J^0$ and, consequently,
$p_1,\dots,p_k\in J^1$. Thus, $p\in CSP_J(v)$. Now (\ref{grafoG})
gives the result.

(\ref{condicionK}) follows directly from (\ref{grafoF}), (\ref{grafoG}) and
(\ref{grafoJ}).
\end{proof}

\section{Exchange Leavitt path algebras}

We will say that a  (not-necessarily unital) ring $R$ is an \emph{exchange
ring} (see \cite{Ara1}) if for every element
$x\in R$ the equivalent conditions in the next lemma are satisfied.

\begin{lem}(\cite[Lemma 1.1]{Ara1})\label{exch def} Let $ R $ be a ring and
let $ R^\prime$ be a unital ring containing $ R $ as a two-sided
ideal. Then the following conditions are equivalent for an element
$ x \in R :$
\begin{enumerate}
\item There exists $e^2=e \in R $ with $e-x \in R^\prime (x-x^2),$
\item there exist $e^2=e \in R  x $ and $ c\in R^\prime $ such that
$(1-e)-c(1-x) \in J(R^\prime),$
\item there exists $e^2=e \in R  x $ such that $R^\prime =R e + R^\prime
(1-x),$
\item there exists $e^2=e \in R x $ such that $1-e \in R^\prime (1-x) ,$
\item there exist $r, s \in R, $ $ e^2=e \in R $ such that $e=rx=s+x-sx.$
\end{enumerate} (Here $J(R^\prime)$ denotes the Jacobson radical of
$R^\prime$.)
\end{lem}

Observe that $R$ being an exchange ring does not depend on the particular
unital ring  where $R$  is embedded as an
ideal (look at condition (v) in the previous Lemma). Other characterizations
of the exchange property for not
necessarily unital rings can be found in \cite{Ara1}.

\begin{rema}\label{ultraexch}
{\rm Since any $K$-matricial algebra is an exchange ring, then so
is any $K$-locally ma\-tricial algebra (apply \cite[Theorem
3.2]{AGS}.)}
\end{rema}

\begin{theor}\label{exch-K}
Let $E$ be a graph. If $L(E)$ is an exchange ring, then $E$ satisfies
condition (K).
\end{theor}
\begin{proof}
The first step will be to show that $E$ satisfies condition (L).
Suppose that there exist a vertex $v$ and a cycle $\alpha$ with
$s(\alpha)=v$ such that $\alpha$ has no exits. Denote by $H$ the
hereditary saturated subset of $E^0$ generated by $\alpha^0$. By
Lemma~\ref{Morita}, $I(H)$ is Morita equivalent to $L(E_H)$. If
$M$ is the graph having only a vertex $w$ and an edge $e$ such
that $r(e)=s(e)=w$, then $L(M)\cong K[x, x^{-1}]$ by \cite[Example
1.4 (ii)]{AA1}. Consider the map $f: L(M)\to L(E_H)$ given by
$f(w)= v$, $f(e)=\alpha$, $f(e^\ast)=\alpha^\ast$. It is well
defined because the relations in $M$ are consistent with those in
$L(E_H)$ (the only non trivial one being $\alpha \alpha^*=v$,
which holds due to the absence of exits for $\alpha$, as in
\cite[p. 12]{AA1}). It is a (non-unital) monomorphism of
$K$-algebras; clearly, $\hbox{Im}f \subseteq v L(E_H) v$. Now, we
prove $vL(E_H)v\subseteq \hbox{Im}f$. To this end, it is enough to
see $vpq^\ast v\in \hbox{Im}f$ for every $p=e^\prime_1\ldots
e^\prime_r$, $q=e_1\dots e_s$, with $e^\prime_1,\dots,e^\prime_r,
e_1,\dots,e_s\in E_H^1$. Reasoning as in \cite[Proof of Theorem
3.11]{AA1} we get that $vpq^\ast v$ has the form: $v, v\alpha^n v$
or $v{(\alpha^\ast)}^m v$, with $m, n \in \N $. Hence our claim
follows.

Apply \cite[Theorem 2.3]{Ara1} to have that $I(H)$ is an exchange
ring; moreover, $L(E_H)$ is an exchange ring by Lemma \ref{Morita}
and \cite[Theorem 2.3]{AGS}, and the same can be said about the
corner $vL(E_H)v$ by \cite[Corollary 1.5]{AGS}. But $vL(E_H)v\cong
L(H)\cong K[x, x^{-1}]$ is not an exchange ring, what leads to a
contradiction.

Now, we will prove that $E$ satisfies condition (K). Suppose on
the contrary that there exists a vertex $v$ and $\alpha=(\alpha_1,
\dots, \alpha_n) \in CSP(v)$, with $\mbox{card}(CSP(v))=1$ (in
fact, $\alpha $ must be a cycle). Consider $A=\{ e\in E^1\mid
e\mbox{ exit of } \alpha\}$, $B=\{ r(e)\mid e\in A\}$, and let $H$
be the hereditary saturated closure of $B$. With a similar
argument to that used in \cite[p. 6]{AA2} we get that $H\cap
\alpha^0=\emptyset$, so that, $H$ is a proper subset of $E^0$.
Then, $\alpha^0 \subseteq (E/H)^0$ and $\{\alpha_1, \dots,
\alpha_n\}\subseteq (E/H)^1$, whence $\alpha$ is a cycle in $E/H$
with no exits.

Since $L(E/H)\cong L(E)/I(H)$ (Lemma \ref{quotient1} (1)),
$L(E/H)$ is an exchange ring \cite[Theorem 2.2]{Ara1} and, by the
previous step, $E/H$ satisfies condition (L), a contradiction.
\end{proof}

Recall that an idempotent $e$ in a ring $R$ is called
\emph{infinite} if $eR$ is isomorphic as a right $R$-module to a
proper direct summand of itself. The ring $R$ is called
\emph{purely infinite} in case every right ideal of $R$ contains
an infinite idempotent.

\begin{prop}\label{K-exch_fin}
If $E$ is a  graph satisfying condition (K) and
$\mathcal{L}(L(E))$ is finite, then $L(E)$ is an exchange ring.
\end{prop}
\begin{proof}

Since $\mathcal{L}(L(E))$ is finite, we can construct an ascending chain of
ideals
$$0=I_0\subseteq I_1\subseteq \cdots \subseteq I_n=L(E)$$
such that, for every $0\leq i\leq n-1$, $I_i$ is maximal among the
ideals of $L(E)$ contained in $I_{i+1}$. Now, let us prove the
result by induction on $n$.

If $n=1$, then $L(E)$ is a simple ring and then $E$ is cofinal by Lemma
\ref{cofinohersat} and \cite[Theorem
3.11]{AA1}. Since $E$ satisfies condition (K), it can occur exactly two
possibilities:
\begin{enumerate}
\item $E$ has no closed simple paths, whence it is acyclic and thus, by
Corollary \ref{ultramatricial}, $L(E)$ is a locally matricial
algebra, and so an exchange ring by Remark \ref{ultraexch}.
\item $E$ has at least one closed simple path, whence $L(E)$ is a purely
infinite simple ring by cofinality,
\cite[Theorem 3.11]{AA1} and \cite[Theorem 11]{AA2}. Thus, $L(E)$ is an
exchange ring by \cite[Corollary 1.2]{Ara2}.
\end{enumerate}
In any case, $L(E)$ turns out to be an exchange ring.

Now, suppose that the result holds for  $k<n$. By Proposition
\ref{iso_ideal} and \cite[Theorem 4.3]{AMFP}, there exist
hereditary saturated sets $H_i$ $(1\leq i\leq n)$ such that:
\begin{enumerate}
\item[(i)] $I_i=I(H_i)$ for every $0\leq i\leq n$; in particular,
$H_i\varsubsetneq H_{i+1}$ for every $0\leq i\leq n-1$.
\item[(ii)] For any $0\leq i\leq n-1$, it does not exist an hereditary
saturated set $T$ such that
$H_i\varsubsetneq T\varsubsetneq H_{i+1}$.
\end{enumerate}
Consider the restriction graph $E_{H_{n-1}}$. By Lemma
\ref{sub_quot_graph}, $E_{H_{n-1}}$ satisfies condition (K), so
that $\mathcal{L}_{\text{gr}}(L(E_{H_{n-1}}))=
\mathcal{L}(L(E_{H_{n-1}}))$ by Proposition \ref{iso_ideal}. If
for each $0\leq i\leq n-1$, $J_i\lhd L(E_{H_{n-1}})$ is the ideal
generated by $H_i$, then the previous remarks imply that
$$0=J_0\subseteq J_1\subseteq \cdots \subseteq J_{n-1}=L(E_{H_{n-1}})$$
where, for every $0\leq i\leq n-2$, $J_i$ is maximal among the
ideals of $L(E)$ contained in $J_{i+1}$; otherwise, since
$\mathcal{L}_{\text{gr}}(L(E_{H_{n-1}}))=
\mathcal{L}(L(E_{H_{n-1}}))$, Lemma \ref{projsather} would
contradict property (ii) satisfied by the set $H_i$. Thus, by
induction hypothesis, $L(E_{H_{n-1}})$ is an exchange ring. Since
$I(H_{n-1})$ is Morita equivalent to $L(E_{H_{n-1}})$ by Lemma
\ref{Morita}, $I(H_{n-1})$ is an exchange ideal by \cite[Theorem
2.3]{AGS}. Now, by Lemma \ref{quotient1} (1),
$L(E)/I(H_{n-1})\cong L(E/H_{n-1})$. Hence, $E/H_{n-1}$ is a graph
satisfying condition (K) by Lemma \ref{sub_quot_graph}, and
$L(E/H_{n-1})$ is simple by construction. Following the same
dichotomy for $E/H_{n-1}$ as in (1) and (2) above, we get that
$L(E/H_{n-1})$ is an exchange ring. Then, by using Lemma
\ref{gen_epi} and \cite[Theorem 3.5]{Ara1}, we conclude that
$L(E)$ is an exchange ring, as desired.
\end{proof}

We would like to thank Gene Abrams for showing that $(4) \Rightarrow (3)$ in
the following theorem is true.

\begin{theor}\label{tres condiciones}
For a graph $E$, the following conditions are equivalent:
\begin{enumerate}
\item $L(E)$ is an exchange ring.
\item $E/H$ satisfies condition (L) for every hereditary saturated subset
$H$ of $E^0$.
\item $E$ satisfies condition (K).
\item $\mathcal{L}_{\text{gr}}(L(E))=\mathcal{L}(L(E))$.
\item $E_H$ and $E/H$ satisfy condition (K) for every hereditary saturated
subset $H$ of $E^0$.
\item $E_H$ and $E/H$ satisfy condition (K) for some hereditary saturated
subset $H$ of $E^0$.
\end{enumerate}
\end{theor}
\begin{proof}$\mbox{ }$

$(1)\Rightarrow (2)$. By Lemma \ref{quotient1} (1), $L(E)/I(H)\cong L(E/H)$.
Then, by \cite[Theorem 2.2]{Ara1},
$L(E/H)$ is an exchange ring. Apply Theorem \ref{exch-K} and Remark
\ref{nova2} (1) to obtain (2).

$(2)\Rightarrow (3)$ is just the first paragraph in the proof of Theorem
\ref{exch-K}.

$(3)\Rightarrow (4)$ is Proposition \ref{iso_ideal}.

$(4)\Rightarrow (3)$. Suppose on the contrary that $E$ does not
satisfy condition (K). Apply $(2)\Rightarrow (3)$ to find a
hereditary saturated subset $H$ of $E^0$ such that $E/H$ does not
verify condition (L), that is, there exists a cycle $p$ in $E/H$
based at $v$ without an exit. Now \cite[Theorem 3.11, pp. 12,
13]{AA1} shows that in this situation we have $v\not\in
J:=I(v+p)$, meaning in particular that the ideal $J$ is not
graded. Now if $H\neq \emptyset$, Lemma \ref{quotient1} shows that
there exists a graded isomorphism $\Phi: L(E)/I(H) \to L(E/H)$ so
that we can lift $\Phi^{-1}(J)$ to an ideal $\mathcal{J}$ of
$L(E)$, which cannot be graded (a quotient of graded ideals is
again graded). If $H=\emptyset$ then clearly $J$ is an ideal of
$L(E/H)=L(E)$ which is not graded. In any case we get a
contradiction.

$(3)\Rightarrow (1)$. We have two different proofs of this fact.
The first one is inspired in the results of \cite{AMFP}, while the
second one follows the style of \cite[Proof of Theorem 4.1]{JP}.
\begin{enumerate}
\item[(i)] By Lemma \ref{finalkey}, there exists a family $\{X_n\}_{n\geq
0}$ of finite subgraphs such that, for every
$n\geq 0$, $X_n$ satisfies condition (K), $E=\bigcup_{n\geq 0}X_n$ and the
natural inclusion maps
$f_n:X_n\hookrightarrow E$ are complete graph homomorphisms (therefore so
are the inclusions
$f_{n,n+1}:X_n\hookrightarrow X_{n+1}$). By \cite[Lemma 2.2]{AMFP}, we have
induced maps
$L(f_{n,n+1}):L(X_n)\rightarrow L(X_{n+1})$ and $L(f_n):L(X_n)\rightarrow
L(E)$ such that $L(E)\cong \varinjlim
(L(X_n), L(f_{n,n+1}))$.

Fix $n\geq 0$. Since $X_n$ satisfies condition (K), by Proposition
\ref{iso_ideal} and \cite[Theorem 4.3]{AMFP},
$\mathcal{L}(L(X_n))$ is isomorphic to the lattice of hereditary saturated
subsets of $X_n^0$. Hence,
$\mathcal{L}(L(X_n))$ is finite. Thus, $L(X_n)$ is an exchange ring by
Proposition \ref{K-exch_fin}. Since $L(E)$ is a
direct limit of exchange rings, it is itself an exchange ring, as desired.

\item[(ii)] Take an element $x\in L(E)$. By \cite[Lemma 1.5]{AA1},
there exist a finite family of vertices $V=\{v_1, \dots, v_m\}$, and a
finite family of edges $W=\{e_1, \dots, e_n\}$,
such that $x$ is in the linear span of $V$ and the set consisting of
expressions $e_{i_1}\dots e_{i_r}e^\ast_{j_1}\dots
e_{j_s}^\ast$, with $e_{i_l}$ edges in $W$. For $T^1=W$ and $T^0= V\cup
\{r(e_i),s(e_i), i=1, \dots, n\}$, let $T$ be
the graph $T=(T^0, T^1, r\vert_{T^1}, s\vert_{T^1})$ and consider $G$, $J$
and $S$ as in Lemma
\ref{conditionKcompletions}. It can be proved, as in \cite[p. 224]{JP}, that
the number of hereditary subsets of $J^0$
is finite. Apply Lemma \ref{conditionKcompletions} (\ref{condicionK}),
Proposition \ref{iso_ideal} and \cite[Theorem
4.3]{AMFP} to obtain that $\mathcal{L}(L(J))$ is finite. Now, Proposition
\ref{K-exch_fin} shows that $L(J)$ is an
exchange ring. By Lemma \ref{quotient1} (1), $L(J)\cong
L(G)/I(\overline{S})$, and by Lemma \ref{gen_epi} and
\cite[Theorem 3.5]{Ara1}, $L(G)$ is an exchange ring. This means (use
condition (5) in Lemma \ref{exch def}) that given
$x$ there exist $e^2=e, r, s\in L(G)\subseteq L(E)$ (Lemmas
\ref{finitecompletions} (\ref{cond1}) and \ref{medkey})
such that $e=rx=s+x-sx$. Whence, $L(E)$ is an exchange ring.
\end{enumerate}

$(3)\Rightarrow (5)$ is Lemma \ref{sub_quot_graph}.

$(5)\Rightarrow (6)$ is a tautology.

$(6)\Rightarrow (1)$. By $(3) \Rightarrow (1)$, $L(E_H)$ and
$L(E/H)$ are exchange rings. Since $L(E_H)$ is Morita equivalent
to $I(H)$ by Lemma \ref{Morita} then $I(H)$ is an exchange ring
because both are idempotent rings and we may apply \cite[Theorem
2.3]{AGS}. By Lemma \ref{quotient1} (1), $L(E/H)\cong L(E)/I(H)$.
Now, $L(E)/I(H)$ and $I(H)$ exchange rings, Lemma \ref{gen_epi}
and \cite[Theorem 3.5]{Ara1} imply that $L(E)$ is an exchange
ring.
\end{proof}



\section{Some special facts}

The following definitions are particular cases  of those appearing
in \cite[Definition 1.3]{DHSz}:

Let $E$ be a graph, and let $\emptyset \ne H\in \mathcal{H}_E$.
Define
$$F_E(H)=\{ \alpha =(\alpha_1, \dots ,\alpha_n)\mid \alpha _i\in
E^1, s(\alpha _1)\in E^0\setminus H, r(\alpha _i)\in E^0\setminus
H \mbox{ for } i<n, r(\alpha _n)\in H\}.$$  Denote by
$\overline{F}_E(H)$ another copy of $F_E(H)$. For $\alpha\in
F_E(H)$, we write $\overline{\alpha}$ to denote a copy of $\alpha$
in $\overline{F}_E(H)$. Then, we define the graph ${}_HE=({}_HE^0,
{}_HE^1, s', r')$ as follows:
\begin{enumerate}
\item ${}_HE^0=({}_HE)^0=H\cup F_E(H)$. \item ${}_HE^1=({}_HE)^1=\{ e\in
E^1\mid s(e)\in H\}\cup \overline{F}_E(H)$.
\item For every $e\in E^1
\mbox{ with } s(e)\in H$, $s'(e)=s(e)$ and $r'(e)=r(e)$. \item For
every $\overline{\alpha}\in \overline{F}_E(H)$,
$s'(\overline{\alpha})=\alpha$ and
$r'(\overline{\alpha})=r(\alpha)$.
\end{enumerate}

\begin{lem}\label{properties}
Let $E$ be a graph, and let $\emptyset \ne H\in \mathcal{H}_E$.
Then:
\begin{enumerate}
\item If $E_H$ satisfies condition (L), then so does ${}_HE$.
\item If $E_H$ satisfies condition (K), then so does ${}_HE$.
\end{enumerate}
\end{lem}
\begin{proof}
Notice that each vertex $\alpha\in {F_E(H)}$ is a source emitting
exactly one edge $\overline{\alpha}\in \overline{F}_E(H)$ which
ends in $H$. Thus, every closed simple path in the graph ${}_HE$
comes from a closed simple path in $E_H$, hence, the result
follows.
\end{proof}

The class of Leavitt path algebras is closed under quotients
(Lemma \ref{quotient1}(1)). A direct consequence of the next
result is that under condition (L), this class is also closed for
ideals.

\begin{lem}\label{isoideal}
{\rm (c.f. \cite[Lemma 1.5]{DHSz})} Let $E$ be a graph, and let
$\emptyset \ne H\in \mathcal{H}_E$. If $E_H$ satisfies condition
(L), then $I(H)\cong L({}_HE)$ (as nonunital rings).
\end{lem}
\begin{proof}
We define a map $\phi: L({}_HE)\rightarrow I(H)$ as follows: (i)
For every $v\in H$, $\phi (v)=v$; (ii) for every $\alpha \in
F_E(H)$, $\phi (\alpha )=\alpha \alpha ^*$; (iii) for every $e\in
E^1 \mbox{ with } s(e)\in H$, $\phi (e)=e$ and $\phi (e^*)=e^*$;
(iv) for every $\overline{\alpha}\in \overline{F}_E(H)$, $\phi
(\overline{\alpha})=\alpha$ and $\phi
(\overline{\alpha}^*)=\alpha^*$.

By definition, it is tedious but straightforward to check that the
images of the relations in $L({}_HE)$ satisfy the relations
defining $L(E)$. Thus, $\phi$ is a well-defined $K$-algebra
morphism.

Since for any $v\in H$, $\phi (v)=v$, to see that $\phi$ is
surjective, by \cite[Lemma 1.5]{AA1}, it is enough to show that
every finite path $\alpha$ of $E$ with $r(\alpha)$ or $s(\alpha)$
in $H$ is in the image of $\phi$. So let $\alpha =(\alpha_1, \dots
,\alpha_n)$ be with $\alpha _i\in E^1$. If $s(\alpha)\in H$, then
$s(\alpha_i)\in H$ for every $i$ because $H$ is hereditary and
thus $\alpha=\phi(\alpha _{1})\cdots \phi (\alpha
_n)=\phi(\alpha)$.

Suppose that $s(\alpha _1)\in E^0\setminus H$ and $r(\alpha _n)\in
H$. Then, there exists $1\leq j\leq n-1$ such that $r(\alpha
_j)\in E^0\setminus H$ and $r(\alpha _{j+1})\in H$. Thus, $\alpha
=(\alpha_1, \dots ,\alpha_{j+1})(\alpha_{j+2}, \dots ,\alpha_n)$,
where $\beta=(\alpha_1, \dots ,\alpha_{j+1})\in F_E(H)$. Hence,
$\alpha= \phi (\overline{\beta})\phi (\alpha _{j+2})\cdots \phi
(\alpha _n)$.

Analogously it can be proved that $\alpha^* \in \mathrm{Im }
\phi$.

Finally, if $0\ne \mbox{Ker}(\phi)$, then
$\mbox{Ker}(\phi)\cap(_HE)^0\ne \emptyset$ by \cite[Proposition
6]{AA2} and Lemma \ref{properties}(2), contradicting the
definition of $\phi$.
\end{proof}

Note that the isomorphism above is not $\mathbb{Z}$-graded because
while $\overline{\alpha}$ has degree $1$ in $_HE$ for every
$\alpha\in F_E(H)$, $\phi(\overline{\alpha})=\alpha$ has not
necessarily degree $1$.

\begin{lem}\label{inducideal}
Let $E$ be a graph satisfying condition (K). Then:
\begin{enumerate}
\item If $J\lhd I \lhd L(E)$, then $J\lhd L(E)$. \item In
particular, if $H\in \mathcal{H}_E$ and $J\lhd I(H)$, then there
exists $X\in \mathcal{H}_E$ such that $X\subset H$ and $J=I(X)$.
\end{enumerate}
\end{lem}
\begin{proof} (1). By Proposition \ref{iso_ideal},
$I=I(H)$ for $H=I\cap E^0\in \mathcal{H}_E$, and by Proposition
\ref{isoideal} (and Remark \ref{nova2}(1)) $I(H)$ is isomorphic to
the Leavitt path algebra $L(_HE)$ and therefore it has a set of
local units. Take $x\in J$ and $z\in L(E)$, then there exits $y\in
I$ such that $x=xy=yx$. Now $zx=(zy)x\in IJ\subseteq J$ and
analogously $xz\in J$.

(2). Again Proposition \ref{iso_ideal} gives that $J=I(X)$ for
$X=J\cap E^0$, and therefore $X=J\cap E^0\subseteq I\cap E^0=H$.
\end{proof}

\begin{prop}\label{noquot}
Let $E$ be a graph satisfying condition (K), let $$X_0=\{v\in E^0\mid
\exists e\ne f\in E^1 \mbox{ with } s(e)=s(f)=v,
r(e)\geq v, r(f)\geq v\},$$ and let $X$ be the hereditary saturated closure
of $X_0$. If $L(E)$ has no unital purely
infinite simple quotients, then neither has $I(X)$.
\end{prop}
\begin{proof}
We will suppose that $X_0\neq \emptyset$, because otherwise there is nothing
to prove.

\underline{Case 1}. We will begin by proving that if $L(E)$ has no
\pis quotients, then $I(X)$ cannot be a \pis ring. Suppose that
this statement is false. By Lemma \ref{isoideal} and Remark
\ref{nova2} (1), $I(X)\cong L({}_XE)$, thus, since $I(X)$ is
unital, ${}_XE$ is a finite graph; in particular, both $X$ and
$F_E(X)$ are finite, and so are
$$X_1=\{v\in E^0\mid v=s(\alpha _i) \mbox{ for some }\alpha
=(\alpha_1, \dots, \alpha_n)\in F_E(X)\}$$ and $Y=X\cup X_1$. We
claim that $K=E^0\setminus Y$ belongs to $\mathcal{H}_E$. Let
$v\in K$, $w\in E^0$, $e\in E^1$ be such that $s(e)=v$ and
$r(e)=w$. We want to prove $w\in K$. Suppose on the contrary that
$w\in Y$. If $w\in X$, then $e\in F_E(X)$ and so $v=s(e)\in
X_1\subseteq Y$, a contradiction, hence $w\in X_1\setminus X$. In
this case there exists a path $\alpha=(\alpha_1, \dots
,\alpha_n)\in F_E(X)$ such that $w=s(\alpha_i)$, for some
$i\in\{1, \dots, n-1\}$. Then $\beta=(e, \alpha_i, \dots
,\alpha_n)\in F_E(X)$ and $v=s(\beta)\in X_1\subseteq Y,$ a
contradiction. This shows that $K$ is hereditary. Now we prove
that it is saturated. Consider $v\in E^0$ and $\emptyset\neq
r(s^{-1}(v))\subseteq K$. Suppose $v\notin K$. Then  $v\in X$ or
$v\in X_1\setminus X$. In the first case, since $X$ is hereditary,
$\emptyset\neq r(s^{-1}(v))\subseteq X$, a contradiction. In the
second one, there exists $\alpha=(\alpha_1, \dots ,\alpha_n)\in
F_E(X)$ such that $v=s(\alpha_i)$ for some $i\in\{1, \dots,
n-1\}$. Then $r(\alpha_i)\in r(s^{-1}(v))\subseteq K$, a
contradiction because $r(\alpha_i)\in Y$, by the definition of
$Y$.

The following step consists of showing that $L(E/K)$, which is isomorphic to
$L(E)/I(K)$ by Lemma \ref{quotient1}, is a
\pis ring. First note that $(E/K)^0=Y$ is finite and therefore $L(E/H)$ is a
unital ring.

Now, since $X$ is finite, $L(E_X)$ is unital. As $L(E_X)$ is
Morita equivalent to the \pis ring $I(X)$ by Lemma \ref{Morita},
$L(E_X)$ is purely infinite simple. By \cite[Proposition 10]{AA2},
$E_X$ is cofinal, satisfies condition (L), and every vertex in
$E_X^0$ connects to a cycle. As $E$ satisfies condition (K), so
does $E/K$ by Lemma \ref{sub_quot_graph}, whence $E/K$ satisfies
condition (L).  Observe that $E/K$ contains at least a cycle;
moreover, since every vertex in $F_E(X)$ connects to a vertex in
$X$, then every vertex in $E/K$ connects to a cycle. Finally,
notice that $E/K$ has no sinks, as otherwise, since any sink would
be in $X$ because $(E/K)^0=X\cup X_1$ and $X_1\setminus X$ clearly
does not have sinks, ${}_XE$ would have a sink, which is not
possible because $L({}_XE)\cong I(X)$ is a \pis ring. Then,
$(E/K)^{\leq \infty}=(E/K)^\infty$. Hence, if $v\in (E/K)^0$ and
 $\alpha$ is in $(E/K)^\infty$, then
$\mbox{card}(\alpha^0)<\infty$ because $Y$ is finite, so $\alpha$
contains a cycle $\beta$  (note that for a cycle $\beta^\prime$
 in $E/K$,  ${\beta^\prime}^0\cap X_0\neq \emptyset$ because $E$
satisfies condition (K), and then by hereditariness
${\beta^\prime}^0\subseteq X$), thus $E/K$ is cofinal. By
\cite[Theorem 11]{AA2} and Lemma \ref{cofinohersat}, $L(E/K)$ is a
\pis ring, a contradiction.

\underline{Case 2}. $I(X)$ has no \pis quotients. Suppose that
$I(X)/J$ is a \pis ring for some ideal $J$ of $I(X)$. By Lemma
\ref{sub_quot_graph}, $E_X$ satisfies condition (K), whence so
does ${}_XE$ by Lemma \ref{properties} (2). Lemma \ref{inducideal}
implies that there exists $H\in \mathcal{H}_E$ such that
$H\subseteq X$ and $J=I(H)$. By Lemma \ref{quotient1}
$L(E)/I(H)\cong L(E/H)$, and by Lemma \ref{sub_quot_graph}, $E/H$
satisfies condition (K). This isomorphism shows that $L(E/H)$ has
no \pis quotients because neither has $L(E)$. If $\Psi$ is the
isomorphism in Lemma \ref{quotient1}, and $Z_0=\Psi(X_0)$, then
$Z=\overline{Z_0}=\Psi (X)$ by Lemma \ref{quotient1}\ (4), and in
particular $I(Z)=\Psi (I(X))$. Thus, by case 1, applied to $E/H,
Z_0$ and $Z$, we get a contradiction.
\end{proof}

The rest of this section is devoted to characterize the primeness
of an ideal of the form $I(H)$, for $H$ hereditary and saturated,
in terms of the so-called maximal tails.

The following definition is a particular case of that of
\cite{BHRS}: Let $E$ be a graph. A nonempty subset $M\subseteq
E^0$ is a \emph{maximal tail} if it satisfies the following
properties:

(MT1) If $v\in E^0$, $w\in M$ and $v\geq w$, then $v\in M$.

(MT2) If $v\in M$ with $s^{-1}(v)\ne \emptyset$, then there exists $e\in
E^1$ with $s(e)=v$ and $r(e)\in M$.

(MT3) For every $v,w\in M$ there exists $y\in M$ such that $v\geq y$ and
$w\geq y$.

\begin{rema}\label{prodideals}
{\rm Let $E$ be a graph. If $J,K\in \mathcal{H}_E$, then
$I(J)I(K)=I(J\cap K)$. To see this, notice that by Remark \ref{gr
id}, $I(J)\cap I(K)=I(J\cap K)$. It is clear that
$I(J)I(K)\subseteq I(J\cap K)$. Since every vertex is an
idempotent, the reverse inclusion is clear.}
\end{rema}

Recall that a graded ideal $I$ of a graded ring $R$ is said to be
\emph{graded prime} if for every pair of graded ideals $J, K$ of
$R$ such that $JK\subseteq I$, it is satisfied that either
$J\subseteq I$ or $K\subseteq I$. The definition of prime ideal is
analogue to the previous one by eliminating the condition of being
graded. It follows by \cite[Proposition II.1.4]{NvO} that for an
ordered group (as it is our case), a graded ideal is graded prime
if and only if it is prime.

\begin{prop}\label{gradedprime}
Let $E$ be a graph, and let $H\in \mathcal{H}_E$. Then, the following are
equivalent:
\begin{enumerate}
\item The ideal $I(H)$ is prime.
\item $M=E^0\setminus H$ is a maximal tail.
\end{enumerate}
\end{prop}
\begin{proof}
$(1)\Rightarrow (2)$. It is not difficult to see that $M$
satisfies (MT1) and (MT2). Suppose that there exist $v,w\in M$
such that no $y\in M$ satisfies:
$$\hbox{
($\ast$)\ $v\geq y$ and $w\geq y$. }$$
 Fix such $v,w$. We will prove
that $\overline{\{v\}}\cap \overline{\{w\}}\cap M=\emptyset$.
Suppose that this is false. Let $m$ be the smallest number such
that $\Lambda_m(T(v))\cap \overline{\{w\}}\cap M\neq \emptyset$
and take $y\in\Lambda_m(T(v))\cap \overline{\{w\}}\cap M$. If
$m>0$, then $s^{-1}(y)\neq \emptyset$ and  $\emptyset \neq
r(s^{-1}(y))\subseteq \Lambda_{m-1}(T(v))\cap \overline{\{w\}}$
because $\overline{\{w\}}$ is hereditary. By the minimality of
$m$, $\Lambda_{m-1}(T(v))\cap \overline{\{w\}}\cap M=\emptyset$,
hence $r(s^{-1}(y))\subseteq M$. Since $M$ is saturated, this
implies $y\in M$, a contradiction. Analogously it can be proved
that $0$ is the smallest number $n$ such that  $T(v)\cap
\Lambda_n(T(w))\cap M\neq \emptyset$, that is,  $T(v)\cap T(w)\cap
M\neq \emptyset $, but this is a contradiction by ($\ast$). Now,
$I(v)I(w)=$ (by Lemma \ref{qualsevol})$
I(\overline{\{v\}})I(\overline{\{w\}})=$ (by Remark
\ref{prodideals})
$I(\overline{\{v\}}\cap\overline{\{w\}})\subseteq$ (as we have
just proved) $I(H)$. By (1), and taking into account Lemma
\ref{qualsevol}), this implies $v\in I(H)$ or $w\in I(H)$, a
contradiction.

 $(2)\Rightarrow (1)$.  Consider two ideals $J_1$ and $J_2$ in
 $L(E)$ such that  $J_1J_2\subseteq I(H)$. By Remark \ref{gr id}
 there exist $H_1, H_2\in \mathcal{H}_E$ such that $J_1=I(H_1)$
 and $J_2=I(H_2)$. By Remark \ref{prodideals}, $H_1\cap H_2\subseteq H$. If
$H_i\varsubsetneq H$ for $i=1,2$, then there exist $v_i\in
H_i\setminus H$ ($i=1,2$). In particular, $v_1,v_2\in M$, so that
there exists $x\in M$ such that $v_i\geq x$ ($i=1,2$). Hence,
$x\in H_1\cap H_2\subseteq H$, which contradicts $x\in M$. Thus,
either $H_1\subseteq H$ or $H_2\subseteq H$, and thus either
$I(H_1)\subseteq I(H)$ or $I(H_2)\subseteq I(H)$, as desired.
\end{proof}

\begin{corol}\label{exchange}
If $E$ is a graph satisfying condition (K), then there is a bijection
between maximal tails and prime ideals. In
particular, if $E$ has no proper maximal tails, then $L(E)$ is simple.
\end{corol}
\begin{proof}
The first statement is a consequence of Proposition
\ref{gradedprime} and Proposition \ref{iso_ideal}. This implies
the second statement because the absence of proper maximal tails
is equivalent to the absence of nonzero prime ideals.
\end{proof}

\section{Stable rank for quasi stable rings}

Let $S$ be any unital ring containing an associative ring $R$ as a
two-sided ideal. The following definitions can be found in
\cite{Vas}. A column vector $b=(b_i)_{i=1}^n$ is called
\emph{$R$-unimodular} if $b_1-1,b_i\in R$ for $i>1$ and there
exist $a_1-1,a_i\in R$ ($i>1$) such that $\sum_{i=1}^n a_ib_i=1$.
The \emph{stable rank} of $R$ (denoted by $sr(R)$) is the least
natural number $m$ for which for any $R$-unimodular vector
$b=(b_i)_{i=1}^{m+1}$ there exist $v_i\in R$ such that the vector
$(b_i+v_ib_{m+1})_{i=1}^m$ is $R$-unimodular. If such a natural
$m$ does not exist we say that the stable rank of $R$ is infinite.

Recall that a ring $R$ is said to be \emph{stable} if $R\cong
M_{\infty }(R)$. In this section, we cover the final step of the
proof of Lemma \ref{ranktwo}. To this end, we need to compute the
stable rank of some rings with local units whose behaviour is
similar to that of stable rings with local units. It is not known
if the property we consider should be equivalent to stability of
the ring.

\begin{lem}\label{quasistable}
Let $R$ be a ring with ascending local unit $\{ p_n \}_{n\geq 1}$. If for
every $n\geq 1$ there exists $m>n$ such that
$p_n\lesssim p_m-p_n$, then $\mbox{sr}(R)\leq 2$.
\end{lem}
\begin{proof}
Fix $S$ a unital ring containing $R$ as two-sided ideal. Let $a_1,
a_2, a_3, b_1, b_2, b_3\in S$ such that $a_1-1, a_2, a_3, b_1-1,
b_2, b_3\in R$, while $a_1b_1+a_2b_2+a_3b_3=1$. By hypothesis,
there exists $n\in \N$ such that $a_1-1, a_2, a_3, b_1-1, b_2,
b_3\in p_nRp_n$. Let $m>n$ such that $p_n\lesssim p_m-p_n$. Then,
there exists $q_n\sim p_n$, $q_n\leq p_m-p_n$. In particular,
$q_np_n=p_nq_n=0$. Now, there exist $u\in p_nRq_n$, $v\in q_nRp_n$
such that $uv=p_n$, $vu=q_n$, $u=p_nu=uq_n$ and $v=q_nv=vp_n$.

Fix $v_1=0$, $v_2=u$, $c_1=b_1$, and $c_2=b_2+vb_3$. Notice that
$(a_1+a_3v_1)-1, c_1-1, (a_2+a_3v_2), c_2 \in R$. Also,
$a_3uvb_3=a_3p_nb_3=a_3b_3$, $a_3ub_2=a_3uq_np_nb_2=0$, and
$a_2vb_3=a_2p_nq_nvb_3=0$. Hence,
$$(a_1+a_3v_1)c_1+(a_2+a_3v_2)c_2=a_1b_1+a_2b_2+a_3b_3=1.$$
Thus, any unimodular $3$-row is reducible, whence the result holds.
\end{proof}

A monoid $M$ is {\it cancellative} if whenever $x+z=y+z$, for
$x,y,z\in M$, then $x=y$. And $M$ is said to be {\it unperforated}
in case for all elements $x,y\in M$ and all positive integers $n$,
we have $nx\leq ny$ implies $x\leq y$.

 {\rm Given an abelian monoid $M$, and an element $x\in M$, we
define $$S(M,x)=\{f:M\rightarrow [0,\infty]\mid \mbox{ additive
map such that } f(x)=1\}.$$ Standard arguments show that, when $M$
is a cancellative monoid, then $S(M,x)$ is nonempty for every
nonzero element $x\in M$.}

\begin{lem}\label{nomesale}
Let $R$ be a nonunital ring with ascending local unit
$\{p_n\}_{n\geq 1}$ such that $V(R)$ is cancellative and
unperforated, and let $S_R =\{ s:V(R)\rightarrow \mathbb{R}^+\mid
\mbox{ morphisms of monoids}\}$. If for every $s\in S_R$,
$\sup_{n\geq 1}\{ s([p_n])\}=\infty$, then for every $n\geq 1$
there exists $m>n$ such that $p_n\lesssim p_m-p_n$.
\end{lem}
\begin{proof}
Fix $n\in\N$, and consider $S_n=S(V(R),2[p_n])$. For every $t\in
S_n$, $\sup_{m\geq 1}t([p_m])=\infty$. Otherwise, there exists
$t\in S_n$ such that $\sup_{m\geq 1}t([p_m])=\alpha \in \R^+$.
Since $\{p_n\}_{n\geq 1}$ is a local unit, we conclude that $t(x)<
\infty $ for every $x\in V(R)$, so that $t\in S_R$, contradicting
the hypothesis. Thus, the maps $\widehat{p_k}:S_n\rightarrow
[0,\infty]$, defined by evaluation, satisfy that the (pointwise)
supremum $\sup_{k\geq 1} \widehat{p_k}=\infty$. Since $S_n$ is
compact, there exists $m> n$ such that $1<\widehat{p_m}$, i.e. for
every $s\in S_n$, $s(2[p_n])<s([p_m])$.

Now, take $t\in S(V(R),[p_m])$. Since $p_n< p_m$, $0\leq
t(2[p_n])=a\leq 2$. If $a=0$, then clearly
$0=t(2[p_n])<t([p_m])=1$. If $a\ne 0$, then $t'(-):=a^{-1}\cdot
t(-)$ belongs to $S_n$, whence $1=t'(2[p_n])<t'([p_m])$ by the
argument above. So, $t(2[p_n])<t([p_m])=1$. Thus, for every $t\in
S(V(R), [p_m])$, we have $t(2[p_n])<t([p_m])=1$. By
\cite[Proposition 3.2]{Ror}, $2p_n\lesssim p_m=p_n+(p_m-p_n)$.
Then, since $V(R)$ is cancellative, we get $p_n\lesssim p_m-p_n$,
as desired.
\end{proof}

\begin{defi}\label{leftinf}
{\rm Let $E$ be a graph. For every $v\in E^0$, we define $L(v)=\{
w\in E^0\mid w\geq v\}$. We say that $v\in E^0$ is {\it left
infinite} if $\mbox{card}(L(v))=\infty$.}
\end{defi}

\begin{defis}\label{graphtrace} {\rm Let $E$ be a graph. A {\it graph
trace} on $E$ is a function $g:E^0\rightarrow \R^+$ such that, for
every $v\in E^0$ with $s^{-1}(v)\ne\emptyset$,
$g(v)=\sum\limits_{s(e)=v}g(r(e))$. We define the {\it norm} of
$g$ to be the (possibly infinite) value $\Vert g\Vert
=\sum\limits_{v\in E^0}g(v)$. We say that $g$ is {\it bounded} if
$\Vert g\Vert <\infty$.}
\end{defis}

\begin{rema}\label{graphstate}
{\rm Let $E$ be a graph, let $E^0=\{v_i\mid i\geq 1\}$, let
$p_n=\sum\limits_{i=1}^{n}v_i$, and let $$S_E=\{
s:V(L(E))\rightarrow \mathbb{R}^+\mid \mbox{ morphisms of monoids}\}.$$ By
\cite[Theorem 2.5]{AMFP}, any element $s\in
S_E$ induces a graph trace by the rule $g_s(v)=s([v])$. Moreover, $g_s$ is
bounded if and only if $\sup_{n\in\N}\{
s([p_n])\}<\infty$.

Conversely, by \cite[Theorem 2.5]{AMFP} and \cite[Lemma 3.3]{AMFP}, if $g$
is a graph trace on $E$, and $v,w\in E^0$
with [$v]=[w]\in V(L(E))$, then $g(v)=g(w)$. So, the rule $s_g([v])=g(v)$ is
well-defined and extends by additivity to
an element $s_g\in S_E$. Certainly, $g$ is bounded if and only if
$\sup_{n\in\N}\{ s_g([p_n])\}<\infty$.}
\end{rema}

Next result in the context of $C^*$-algebras is \cite[Lemma
3.8]{Tomf}. Here, we follow a different approach to prove it.

\begin{lem}\label{compquot}
Let $E$ be a graph, let $H\in \mathcal{H}_E$, and let $\pi :L(E)\rightarrow
L(E)/I(H)$ be the natural projection map.
If $e\in L(E)$ is an idempotent, $W\subseteq E^0\setminus H$ is a finite
set, and $\pi (e)\lesssim \sum_{w\in W}\pi
(w)$ in $L(E)/I(H)$, then there exists a finite set $X\subseteq H$ such that
$e\lesssim \sum_{w\in W}w+\sum_{x\in X}x$.
\end{lem}
\begin{proof}
By \cite[Theorem 2.5 and Lemma 5.6]{AMFP}, $V(L(E))/V(I(H))\cong
V(L(E/H))$. Thus, $[\pi (e)]\leq \sum_{w\in W}[\pi (w)]\in
V(L(E/H))$ implies that there exist $a,b\in V(I(H))$ such that
$[e]+a\leq \sum_{w\in W}[w]+b\in V(L(E))$. Since $V(I(H))=\langle
[v]\mid v\in H\rangle$, there exists a finite set $X\subseteq H$
such that $b=\sum_{x\in X}x$. Then, $[e]\leq \sum_{w\in
W}[w]+\sum_{x\in X}[x]$, as desired.
\end{proof}

\begin{prop}\label{tomforde}{\rm (c.f. \cite[Theorem 3.2]{Tomf})}
Let $E$ be a graph. If every vertex of $E$ lying on a closed
simple path is left infinite and $E$ has no nonzero bounded graph
traces, then for every finite set $V\subseteq E^0$ there exists a
finite set $W\subseteq E^0$ with $V\cap W=\emptyset$ and
$\sum_{v\in V}v\lesssim \sum _{w\in W}w$.
\end{prop}
\begin{proof}
The proof of this result corresponds to $(d)\Rightarrow
(e)\Rightarrow (f)$ of \cite[Theorem 3.2]{Tomf}, with suitable
adaptation of the arguments except for the Case II in
$(d)\Rightarrow (e)$, in which the way to prove the following
statement is different: If $F\subseteq E^0$ is a finite set, and
$n=\max\{i\in\N\mid w_i\in F\}$, there exists $m>n$ such that
$p_n\lesssim p_m-p_n$.

Suppose then $v\not\in\overline{H}$. List the vertices of
$E/\overline{H}=\{w_i\mid i\geq 1\}$, in such a way that $w_1=v$.
Let $\pi :L(E)\rightarrow L(E)/I(H)$ be the natural projection
map. For every $n\geq 1$, set $p_n=\sum_{i=1}^{n}\pi (w_i)$.
Clearly, $\{p_n\}_{n\geq 1}$ is an ascending local unit for
$L(E/\overline{H})$. Since every vertex on a closed simple path is
left infinite, no vertex on $E/\overline{H}$ lies on a closed
simple path. Thus, $E/\overline{H}$ is acyclic, whence
$L(E/\overline{H})$ is locally matricial by Corollary
\ref{ultramatricial}. In particular, $V(L(E/\overline{H}))$ is
cancellative and unperforated. Moreover, since $E$ has no nonzero
bounded graph traces, neither has $E/\overline{H}$. Otherwise, by
Remark \ref{graphstate}, there exists a monoid morphism $s:
V(L(E/\overline{H}))\rightarrow \R^+$ with $\sup_{n\in\N}\{
s([p_n])\}<\infty$. Hence, $s$ induces a monoid morphism
$s\circ\pi :V(L(E))\rightarrow \R^+$ such that $\sum_{v\in
E^0}(s\circ\pi)([v])=\sum_{v\in
E^0\setminus\overline{H}}s([v])<\infty$, consequently there exists
a bounded graph trace on $E$, contradicting the assumption. By
Remark \ref{graphstate} and Lemma \ref{nomesale}, for every $n\geq
1$ there exists $m>n$ such that $p_n\lesssim p_m-p_n$.
\end{proof}

\begin{corol}\label{elrizo}
Let $E$ be a graph. If every vertex of $E$ lying on a closed simple path is
left infinite and $E$ has no nonzero
bounded graph traces, then $\mbox{sr}(L(E))\leq 2$.
\end{corol}
\begin{proof}
Let $E^0=\{v_i\mid i\geq 1\}$, and for each $n\in \N$ consider
$p_n=\sum\limits_{i=1}^{n}v_i$. Then, $\{p_n\}_{n\geq 1}$ is an
ascending local unit for $L(E)$. Fix $n\geq 1$ and set $V=\{v_1,
\dots ,v_n\}$. By Proposition \ref{tomforde}, there exists a
finite subset  $W\subseteq E^0$ such that $V\cap W=\emptyset$ and
$p_n=\sum_{v\in V}v\lesssim \sum _{w\in W}w$. If $m$ is the
largest subindex of $w\in W$, notice that $m> n$ and that $\sum
_{w\in W}w\leq p_m-p_n$. Hence, the result holds because $L(E)$
satisfies the hypotheses of Lemma \ref{quasistable}.
\end{proof}

\section{Stable rank for exchange Leavitt path algebras}

In this section, we characterize the stable rank of exchange
Leavitt path algebras in terms of intrinsic properties of the
graph.

\begin{lem}\label{rankone}
Let $E$ be an acyclic graph. Then, the stable rank of $L(E)$ is $1$.
\end{lem}
\begin{proof} If $E$ is finite, then $L(E)$ is a $K$-matricial algebra by
Corollary \ref{matricial}, whence $\mbox{sr}(L(E))=1$. Now suppose
that $E$ is infinite. By Corollary \ref{ultramatricial}, there
exists a family $\{X_n\}_{n\geq 0}$ of finite subgraphs of $E$
such that $L(E)\cong \varinjlim L(X_n)$. By the definitions of
direct limit and stable rank,
$$(\ast)\qquad \mbox{sr}(L(E))\leq \liminf\limits_{n\rightarrow \infty}
\mbox{sr}(L(X_n)).$$

If $E$ is acyclic, then so are the $X_n$'s, whence $\mbox{sr}(L(E))=1$ by
the result above and $(\ast)$.
\end{proof}

\begin{lem}\label{pisuquot}
Let $E$ be a graph satisfying condition (K). Then, $L(E)$ has a
unital purely infinite simple quotient if and only if there exists
$H\in \mathcal{H}_E$ such that the quotient graph $E/H$ is
nonempty, finite, cofinal and contains no sinks.
\end{lem}
\begin{proof}
First, suppose that $J$ is an ideal of $L(E)$ such that $L(E)/J$
is a unital purely infinite simple ring. By Proposition
\ref{iso_ideal}, there exists $H\in \mathcal{H}_E$ such that
$J=I(H)$. By Lemma \ref{quotient1} (1), $L(E)/J\cong L(E/H)$.
Moreover, $E/H$ satisfies condition (K) by Lemma
\ref{sub_quot_graph}. Hence, since $L(E/H)$ is unital, $E/H$ is
finite. Since $L(E)/J$ is purely infinite simple, $E/H$ is cofinal
and every vertex connects to a closed simple path by \cite[Theorem
11]{AA2}, whence $E/H$ has no sinks.

Conversely, suppose that there exists $H\in \mathcal{H}_E$ such
that the quotient graph $E/H$ is nonempty, finite, cofinal and
contains no sinks. Thus, it contains a closed simple path and
every vertex connects to a closed simple path. Then, since $E/H$
satisfies condition (K) by Lemma \ref{sub_quot_graph}, $L(E/H)$ is
unital, purely infinite and simple by \cite[Theorem 11]{AA2} and
\ref{cofinohersat}. By Lemma \ref{quotient1} (1), $L(E)/M\cong
L(E/H)$ and the proof is complete.
\end{proof}

\begin{corol}\label{rankinf}
Let $E$ be a graph satisfying condition (K). If there exists $H\in
\mathcal{H}_E$ such that the quotient graph $E/H$ is
nonempty, finite, cofinal and contains no sinks, then the stable rank of
$L(E)$ is $\infty$.
\end{corol}
\begin{proof}
By Lemma \ref{pisuquot}, there exists a maximal ideal $M\lhd L(E)$
such that $L(E)/M$ is a unital purely infinite simple ring. Thus,
$\mbox{sr}(L(E)/M)=\infty$ (see \cite{APP}). Since
$\mbox{sr}(L(E)/M)\leq \mbox{sr}(L(E))$ (see \cite[Theorem
4]{Vas}), we conclude that $\mbox{sr}(L(E))=\infty$.
\end{proof}

The proof of the following result closely follows  that of
\cite[Lemma 3.2]{DHSz}.

\begin{lem}\label{ranktwo}
Let $E$ be a non acyclic graph satisfying condition (K). If $L(E)$
does not have any unital purely infinite simple quotient, then
there exists a graded ideal $J\lhd L(E)$ with $\mbox{sr}(J)=2$
such that $L(E)/J$ is a locally matricial $K$-algebra.
\end{lem}
\begin{proof}
Let $$X_0=\{v\in E^0\mid \exists e\ne f\in E^1 \mbox{ with }
s(e)=s(f)=v, r(e)\geq v, r(f)\geq v\},$$ and let $X$ be the
hereditary saturated closure of $X_0$. Consider $J=I(X)$, and
notice that $L(E)/J\cong L(E/X)$ by Lemma \ref{quotient1} (1).
Moreover, since $E$ satisfies condition (K), then so does $E/X$ by
Lemma \ref{sub_quot_graph}. If there is a closed simple path
$\alpha $ in $E/X$, then every $v\in \alpha ^0$  satisfy
$\mbox{card}(CSP_{E/X}(v))\geq 2$, therefore, there exists a
vertex $v_0\in \alpha^0\cap X_0\subseteq X$, contradicting the
assumption. So, $E/X$ contains no closed simple paths, whence it
is an acyclic graph, and thus $L(E)/J$ is locally matricial by
Corollary \ref{ultramatricial}.

Now, by  Remark \ref{nova2} (1), Lemma \ref{sub_quot_graph} and
Lemma \ref{isoideal}, $J\cong L({}_XE)$. We will show that every
vertex lying in a closed simple path of ${}_XE$ is left infinite,
and that ${}_XE$ has no nonzero bounded graph traces, as a way of
contradiction.

Suppose that there exists a closed simple path $\alpha$ in ${}_XE$
such that the set $Y$ of vertices of ${}_XE$ connecting to the
vertices of $\alpha ^0$ is finite. It is not difficult to see that
$\alpha^0\cup Y$ is a maximal tail in ${}_XE$. Let $M$ be a
maximal tail of the smallest cardinal contained in $\alpha^0\cup
Y$. Observe that $M\cap X_0\ne \emptyset$; otherwise $X\setminus
M$, which is a hereditary saturated proper subset of $X$, would
contain $X_0$, which is impossible. Denote by $\widetilde{M}$ the
quotient graph of ${}_XE$ by the hereditary saturated set
$H={}_XE^0\setminus M$, i.e. $\widetilde{M}={}_XE/H$. Then, since
$M$ is finite, $L(\widetilde{M})$ is a unital ring. As $E$
satisfies condition (K), so does ${}_XE$ (by Lemma
\ref{sub_quot_graph} and Lemma \ref{properties} (2)) and thus
$\widetilde{M}$ (by Lemma \ref{sub_quot_graph} again). Then, since
$M$ does not contain smaller maximal tails, $L(\widetilde{M})$ is
simple by Corollary \ref{exchange}. As $M\cap X_0\ne \emptyset$,
$\widetilde{M}$ is non acyclic. Thus, $L(\widetilde{M})\cong
L({}_XE)/I(H)$ (by Lemma \ref{quotient1} (1)) is a unital purely
infinite simple ring by \cite[Theorem 3.11]{AA1}, Lemma
\ref{cofinohersat} and \cite[Theorem 11]{AA2}. By Proposition
\ref{noquot} $L(E)$ has a unital purely infinite simple quotient,
contradicting the hypothesis. Hence, every vertex lying in a
closed simple path in ${}_XE$ is left infinite.

Now, suppose that there exists a nonzero bounded graph trace $g$
on ${}_XE$. By Remark \ref{graphstate},
$s_g:V(L({}_XE))\rightarrow \R^+$ is a nonzero morphism such that
$\sum_{v\in {}_XE^0}s_g([v])<\infty$. But for any $v\in X_0$ we
have $2s_g([v])\leq s_g([v])$, so that $g(v)=0$. Hence,
$X_0\subseteq \{w\in {}_XE^0\mid g(w)=0\}$, which is a hereditary
saturated subset of ${}_XE$ by \cite[Lemma 3.7]{Tomf}. Thus, since
${}_XE=\overline{X_0}^{({}_XE)}$, we conclude that $g\equiv 0$,
contradicting the assumption. Hence, there exist no nonzero
bounded graph traces on ${}_XE$.

Thus, $\mbox{sr}(J)=\mbox{sr}(L({}_XE))\leq 2$ by Corollary
\ref{elrizo}. Since every vertex in $X_0$ is properly infinite as
an idempotent of $L({}_XE)$, $\mbox{sr}(L({}_XE))\ne 1$, so that
$\mbox{sr}(J)= 2$, as desired.
\end{proof}

\begin{corol}\label{lodicho}
Let $E$ be a non acyclic graph satisfying condition (K). If $L(E)$ does not
have any unital purely infinite simple
quotient, then $\mbox{sr}(L(E))=2$.
\end{corol}
\begin{proof}
Consider $J$ the graded ideal obtained in the previous Lemma. By
\cite[Theorem 4]{Vas},
$$2=\max\{ \mbox{sr}(J), \mbox{sr}(L(E)/J)\}\leq
\mbox{sr}(L(E))\leq \max\{ \mbox{sr}(J),
\mbox{sr}(L(E)/J)+1\}=2.$$ Then, $\mbox{sr}(L(E))=2$, as desired.
\end{proof}

\begin{theor}\label{trychotomy}
Let $E$ be a graph satisfying condition (K). Then, the values of the stable
rank of $L(E)$ are:
\begin{enumerate}
\item $\mbox{sr}(L(E))=1$ if $E$ is acyclic. \item
$\mbox{sr}(L(E))=\infty$ if there exists $H\in \mathcal{H}_E$ such that the
quotient graph $E/H$ is nonempty, finite,
cofinal and contains no sinks. \item $\mbox{sr}(L(E))=2$ otherwise.
\end{enumerate}
\end{theor}
\begin{proof}
Statement (1) holds by Lemma \ref{rankone}, and statement (2) by
Corollary \ref{rankinf}. If $E$ is non acyclic and it does not
exist $H\in \mathcal{H}_E$ such that the quotient graph $E/H$ is
nonempty, finite, cofinal and contains no sinks, then $L(E)$ does
not have any unital purely infinite simple quotient by Lemma
\ref{pisuquot}. Hence, statement (3) holds by Corollary
\ref{lodicho}.
\end{proof}

\section*{Acknowledgments}

Part of this work was done during a visit of the second author to the
Department of Mathematics of the University of
Colorado at Colorado Springs (U.S.A.) and to the Centre de Recerca
Matem\`atica (U.A.B., Spain). The second author want
to thanks both host centers for their warm hospitality, and in particular to
Gene Abrams and Pere Ara for valuable
discussions on the topics of
this paper.

\end{document}